\documentclass[11pt,a4paper]{article}

\usepackage{color}
\usepackage{graphics,epsfig}
\usepackage{amsfonts,amssymb,amsmath,amsthm}

\newtheorem{lem}{Lemma}[section]
\newtheorem{cor}[lem]{Corollary}
\newtheorem{thm}[lem]{Theorem}
\newtheorem{prop}[lem]{Proposition}
\theoremstyle{definition}
\newtheorem{defi}[lem]{Definition}
\theoremstyle{remark}
\newtheorem{rem}[lem]{Remark}

\numberwithin{equation}{section}

\newcommand{\ep}{\varepsilon}
\newcommand{\ue}{u^\ep}
\newcommand{\ve}{v^\ep}
\newcommand{\fe}{f_\ep}
\newcommand{\de}{\delta}
\newcommand{\U}{U_0}

\newcommand{\VVV}{U_{1zz}}
\newcommand{\hm}{\alpha_-(\delta)}
\newcommand{\h}{a(\delta)}
\newcommand{\hp}{\alpha_+(\delta)}
\newcommand{\m}{\mu(\de)}
\newcommand{\BR}{{\bf R}}
\newcommand{\n}{\nabla }
\newcommand{\dudn}{\displaystyle{\frac{\partial u}{\partial\nu}}}
\newcommand{\dvdn}{\displaystyle{\frac{\partial v}{\partial\nu}}}

\newcommand{\dom}{\partial\Omega}
\newcommand{\om}{\Omega}
\newcommand{\ombar}{\overline{\Omega}}
\newcommand{\edeux}{\displaystyle{\frac{1}{\ep ^2}}}
\newcommand{\R}{\mathbb{R}}

\newcommand{\vsp}{\vspace{8pt}}
\newcommand{\ga}{\gamma}

\newcommand{\di}{\displaystyle}

\newcommand{\f}{f_{\delta}}
\newcommand{\am}{\alpha _ -}
\newcommand{\ap}{\alpha _ +}
\newcommand{\EB}{e^{-\beta t/\ep ^2}}
\newcommand{\Pe}{(P^{\;\!\ep})}
\newcommand{\V}{U_1}

\newcommand{\emuttchemo}{e^{\mu (\ep  G) t/\ep ^2}}
\newcommand{\fde}{f_\de}
\newcommand{\hde}{h_\de}
\newcommand {\Q}{Q_T}
\newcommand{\Y}{Y(\tau,\xi;\de)}
\newcommand{\Yxi}{Y_\xi(\tau,\xi;\de)}
\newcommand{\Yxixi}{Y_{\xi\xi}(\tau,\xi;\de)}

\newcommand{\A}{A(\tau,\xi;\de)}
\newcommand{\dvvdnchemo}{\displaystyle{\frac{\partial {v^0}}{\partial\nu}}}
\newcommand{\gam}{\Gamma}
\newcommand{\x}{\xi}

\title{\bf THE SINGULAR LIMIT OF A CHEMOTAXIS-GROWTH SYSTEM WITH GENERAL INITIAL DATA}
\author{ }
\date{}

\begin{document}
\maketitle \vspace{-20 mm}
\begin{center}
{MATTHIEU ALFARO}\\[1ex]
Analyse Num\'erique et EDP, Universit\'e de Paris Sud,\\
91405 Orsay Cedex, France. \\[2ex]
\end{center}

\vspace{15pt}

\begin{abstract}
We study the singular limit of a system of partial differential
equations which is a model for an aggregation of amoebae subjected
to three effects: diffusion, growth and chemotaxis. The limit
problem involves motion by mean curvature together with a nonlocal
drift term. We consider rather general initial data. We prove a
generation of interface property and study the motion of
interface. We also obtain an optimal estimate of the thickness and
the location of the transition layer that develops.
\end{abstract}

\section{Introduction}\label{intro-chemo}

Let us start by a short description of life-cycles of the cellular
slime molds ({\it amoebae}). The cells feed and divide until
exhaustion of food supply. Then, the amoebae aggregate to form a
multicellular assembly called a slug. It migrates to a new
location, then forms into a fruiting body, consisting of a stalk
formed from dead amoebae and spores on the top (fruiting bodies
that are visible to the naked eye are often referred to as
mushrooms). Under suitable conditions of moisture, temperature,
spores release new amoebae. The cycle then repeats itself.

It is known that the aggregation stage is mediated by {\it
chemotaxis}, i.e. the tendency of biological individuals to direct
their movements according to certain chemicals in their
environment. The chemotactant ({\it acrasin}) is produced by the
amoebae themselves and degraded by an extracellular enzyme ({\it
acrasinase}). For more details on the biological background, we
refer to \cite{KS}, \cite{Nan} or \cite{FL}.

So the amoebae have a random motion analogous to diffusion coupled
with an oriented chemotactic motion in the direction of a positive
gradient of acrasin. In 1970, Keller and Segel \cite{KS} proposed
the following system as a model to describe such movements leading
to slime mold aggregation:
\[(KS) \quad\begin{cases}
u_t&=\n\cdot(D_2\n u)-\n\cdot(D_1\n v),\\
v_t&=D_v \Delta v+f(v)u-k(v)v,
\end{cases}\]
inside a closed region $\om$. Here, $u$, respectively $v$, denotes
the concentration of amoebae, respectively of acrasin; $f(v)$ is
the production rate of acrasin, and $k(v)$ the degradation rate of
acrasin (due to acrasinase); $D_2=D_2(u,v)$, respectively
$D_1=D_1(u,v)$, measures the vigor of the random motion of the
amoebae, respectively the strength of the influence of the acrasin
gradient on the flow of amoebae; $D_v$ is a positive and constant
diffusion coefficient. The problem is completed by initial data
$u_0$ and $v_0$ and, assuming that there is now flow of the
amoebae or the acrasin across the boundary $\dom$, by homogeneous
Neumann boundary conditions
$$
\n u\cdot \nu=\n v\cdot \nu=0 \quad\ \textrm{on} \ \; \dom\times
(0,+\infty),
$$
$\nu$ being the unit outward normal to $\dom$.

An often used simplified model is obtained as follows. By some
receptor mechanism, cells do not measure the gradient of $v$ but
of some $\chi(v)$, with a sensitive function $\chi$ satisfying
$\chi '>0$, so that $D_1(u,v)=u\chi '(v)$. By taking $D_2$, $f$
and $k$ as constant functions and using some rescaling arguments,
the system reduces to
\[(KS\, ') \quad\begin{cases}
\;\;u_t&=d_u\Delta u-\n\cdot(u\n \chi(v)),\\
\tau v_t&=d_v\Delta v+u-\gamma v,
\end{cases}\]
with $d_u$, $d_v$, $\tau$ and $\gamma$ some positive constants.

Many analyses of the Keller-Segel model for the aggregation
process were proposed. Chemotaxis having some features of
``negative diffusion", Nanjundiah \cite{Nan} suggests that the
whole population concentrates in a single point; we refer to this
phenomenon as the {\it chemotactic collapse}. In mathematical
terms, this amounts to blow up in finite time. As a matter of
fact, it turns out that the possibility of collapse depends upon
the space dimension. In particular it never happens in the
one-dimensional case whereas in two space dimensions, assuming
radially symmetric situations, it only occurs if the total amoebae
number is sufficiently large. The problem of global existence and
blow up of solutions has been intensively studied; we refer in
particular to \cite{CP}, \cite{Sch}, \cite{LNT}, \cite{JL},
\cite{N}, \cite{HV1}, \cite{HV2}.

In a different framework, Mimura and Tsujikawa \cite{MT}, consider
aggregating pattern-dy\-na\-mics arising in the following
chemotaxis model with growth:
\[(MT^\ep) \quad\begin{cases}
\;\; u_t&=\ep ^2 \Delta u-\ep \n \cdot (u\n \chi (v))+f(u),\\
\tau v_t&=\Delta v+u-\gamma v,
\end{cases}\]
where $\ep >0$ is a small parameter. The function $f$ is cubic, 0
and 1 being its stable zeros, and satisfies $\int _0 ^1 f>0$. In
this model, the population is subjected to three effects:
diffusion, growth and chemotaxis. The diffusion rate and the
chemotactic rate are both very small compared with the growth
rate. They observe that, in a first stage, internal layers
--- which describe the boundaries of aggregating regions ---
develop; in a second stage, the motion of the aggregating regions
--- which can be described by that of internal layers --- takes
place. The balance of the three effects (diffusion, growth and
chemotaxis) makes the aggregation mechanism possible. Taking the
limit $\ep\to 0$, they formally derive the equation for the motion
of the limit interface and study the stability of radially
symmetric equilibrium solutions.

The purpose of this paper is to extend some of the results
obtained by Bonami, Hilhorst, Logak and Mimura \cite{BHLM3} about
the singular limit of a variant of system $(MT^\ep)$, where the
second equation is elliptic ($\tau=0$):
\[(P^\ep) \quad\begin{cases}
u_t=\Delta u-\n \cdot (u\n \chi(v))+\edeux \fe (u) &\textrm{ in }
\Omega \times
(0,+\infty),\vspace{3pt}\\
0=\Delta v +u - \gamma v &\textrm{ in } \Omega \times
(0,+\infty),\vspace{3pt}\\
\dudn =\dvdn= 0 &\textrm{ on }\partial \Omega \times (0,+\infty),\vspace{3pt} \\
u(x,0)=u_0(x)   &\textrm{ in } \Omega,
\end{cases}\]
where $\om$ is a smooth bounded domain in $\R^N$ ($N\geq 2$),
$\nu$ is the Euclidian unit normal vector exterior to $\partial
\om$. We assume that $\gamma$ is a positive constant and that the
nonlinearity$\fe$ is given by
\begin{equation}\label{nonlinearity-chemo}
\begin{array}{ll}
\fe(u)&=u(1-u)(u-\di{\frac 1 2})+\ep \alpha u(1-u)\vsp\\
&=:f(u)+\ep g(u),
\end{array}
\end{equation}
with $\alpha>0$. The role of the function $g$ is to break the
balance of the two stable zeros slightly. The sensitive function
$\chi$ is smooth and satisfies $\chi'(v)>0$ for $v>0$.

We also assume that the initial datum satisfies $u_0 \in
C^2(\ombar)$ and $u_0 \geq 0$. Throughout the present paper, we
fix a constant $C_0>1$ that satisfies
\begin{equation}\label{int1-chemo}
\|u_0\|_{C^0(\ombar)}+\|\n u_0\|_{C^0(\ombar)}+\| \Delta
u_0\|_{C^0(\ombar)}\leq C_0.
\end{equation}
Furthermore we define the ``initial interface" $\Gamma _0$ by
\[\gam _0:=\{x\in\om, \; u_0(x)=1/2\}.\]
We suppose that $\gam _0$ is a $C^{2+\vartheta}$ hypersurface
without boundary, for a $\vartheta \in (0,1)$, such that, $n$
being the Euclidian unit normal vector exterior to $\gam _0$,
\begin{equation}\label{dalltint-chemo}
\Gamma_0 \subset\subset \Omega \quad \mbox { and } \quad \n
u_0(x)\cdot n(x) \neq 0\quad\text{if $x\in\gam _0,$}
\end{equation}
\begin{equation}\label{initial-data-chemo}
u_0>1/2 \quad\ \text {in}\; \  \om ^{(1)} _0,\quad u_0<1/2 \quad\
\text {in}\; \ \om ^{(0)} _0 ,
\end{equation}
where $\om ^{(1)} _0$ denotes the region enclosed by $\Gamma _0$
and $\om ^{(0)} _0$ the region enclosed between $\partial \om$ and
$\Gamma _0$.

The existence of a unique smooth solution to Problem $\Pe$ is
proved in \cite{BHLM3}, Lemma 4.2:

\begin{lem}\label{bonami-chemo}
There exists $\ep _0 >0$ such that, for all $\ep \in (0,\ep _0)$,
there exists a unique solution $(\ue,\ve)$ to Problem $\Pe$ on
$\om \times [0,+\infty)$, with $0\leq \ue \leq C_0$ on $\Q$.
\end{lem}

To study the interfacial behavior associated with this model, it
is useful to consider a formal asymptotic limit of Problem $\Pe$
as $\ep\rightarrow 0$. Then the limit solution $u^0 (x,t)$ will be
a step function taking the value $1$ on one side of the interface,
and $0$ on the other side. This sharp interface, which we will
denote by $\Gamma _t$, obeys a law of motion, which can be
obtained by formal analysis (see Section
\ref{preliminaries-chemo}):
\[ (P^0)\quad\begin{cases}
V_{n}=-(N-1)\kappa + \di{\frac{\partial {\chi(v^0)}}{\partial n}}
+ \sqrt 2 \alpha
&\textrm{ on } \gam _t,\\
\vspace{3pt}
\gam _t\big|_{t=0}=\gam _0\\
\vspace{3pt} -\Delta v^0+\gamma v^0=u^0 &\textrm{ in } \Omega \times (0,T],\\
\dvvdnchemo= 0  &\textrm{ on }\partial \Omega \times (0,T], \\
\end{cases} \]
where $V_n$ is the normal velocity of $\Gamma _t$ in the exterior
direction, $\kappa$ the mean curvature at each point of
$\Gamma_t$. We set $Q_T:=\om \times [0,T]$ and for each $t\in
[0,T]$, we define $\om_t ^{(1)}$ as the region enclosed by the
hypersurface $\gam _t$ and $\om_t ^{(0)}$ as the region enclosed
between $\partial \om$ and $\gam _t$. The step function $u^0$ is
determined straightforwardly from $\Gamma _t$ by
\begin{equation}\label{u-chemo}
u^0(x,t)=\begin{cases} 1 &\text{in } \om_t ^{(1)}\\
0 &\text{in } \om_t ^{(0)} \end{cases} \quad\text{for } t\in[0,T].
\end{equation}
By a contraction fixed-point argument in suitable H\"o\-lder
spaces, the well-posedness, locally in time, of the free boundary
Problem $(P^0)$ is proved in \cite{BHLM3}, Theorem 2.1:
\begin{lem}\label{bonami-freepb-chemo}
There exists a time $T>0$ such that $(P^0)$ has a unique solution
$(v^0,\Gamma)$ on $[0,T]$, with
$$
\Gamma=\bigcup _{0\leq t \leq T} (\Gamma_t\times\{t\}) \in
C^{2+\vartheta,\frac{2+\vartheta}{2}},
$$
and $v^0 |_ \Gamma  \in C^{2+\vartheta,\frac{2+\vartheta}{2}}$.
\end{lem}

Bonami, Hilhorst, Logak and Mimura \cite{BHLM3} have proved a
motion of interface property; more precisely, for some prepared
initial data,  they show that $(\ue,\ve)$ converges to $(u^0,v^0)$
as $\ep \rightarrow 0$, on the interval $(0,T)$. So the evolution
of $\Gamma _t$ determines the aggregating patterns of the
individuals. Here we consider the case of arbitrary initial data.
Our first main result, Theorem \ref{width-chemo}, describes the
profile of the solution after a very short initial period. It
asserts that, given a virtually arbitrary initial datum $u_0$, the
solution $\ue$ quickly becomes close to $1$ or $0$, except in a
small neighborhood of the initial interface $\Gamma _0$, creating
a steep transition layer around $\Gamma _0$ ({\it generation of
interface}). The time needed to develop such a transition layer,
which we will denote by $t ^\ep$, is of order $\ep^2|\ln\ep|$. The
theorem then states that the solution $\ue$ remains close to the
step function $u^0$ on the time interval $[t^\ep,T]$ ({\it motion
of interface}). Moreover, as is clear from the estimates in the
theorem, the ``thickness" of the transition layer is of order
$\ep$.

\begin{thm}[Generation and motion of interface]\label{width-chemo}
Let $\eta \in(0,1/4)$ be arbitrary and set
$$
\mu=f'(1/2)=1/4.
$$
Then there exist positive constants $\ep _0 $ and $C$ such that,
for all $\,\ep \in (0,\ep _0)$, all $\,t^\ep \leq t \leq T$, where
$t^\ep=\mu ^{-1} \ep ^2 |\ln \ep|$, we have
\begin{equation}\label{resultat-chemo}
\ue(x,t) \in \begin{cases} {}[-\eta,1+\eta]&\quad\text{if}\quad
x\in\mathcal N_{C\ep}(\Gamma_t)\vsp\\
{}[-\eta,\eta]&\quad\text{if}\quad
x\in\om_t^{(0)}\setminus\mathcal N_{C\ep}(\Gamma
_t)\vsp\\
{}[1-\eta,1+\eta]&\quad\text{if}\quad x\in\om
_t^{(1)}\setminus\mathcal N_{C\ep}(\Gamma _t),
\end{cases}
\end{equation}
where $\mathcal N _r(\Gamma _t):=\{x\in \om,dist(x,\Gamma _t)<r\}$
denotes the $r$-neighborhood of $\Gamma _t$.
\end{thm}

\begin{cor}[Convergence]\label{total-chemo}
As $\ep\to 0$, the solution $(\ue,\ve)$ converges to $(u^0,v^0)$
everywhere in $\bigcup _{0<t\leq T}(\om^{(0\ \text{or}\
1)}_t\times\{ t\})$.
\end{cor}

\vskip 8 pt The next theorem deals with the relation between the
set $\Gamma ^\ep _t:=\{x \in \om , u^\ep(x,t)=1/2\}$ and the
solution $\Gamma _t$ of Problem $(P ^0)$.

\begin{thm}[Error estimate]\label{error-chemo}
There exists $C>0$ such that
\begin{equation}\label{thi-2-chemo}
\Gamma _t ^\ep \subset \mathcal N _{C \ep} (\Gamma _t)\quad\ \text
{for}\; \  0 \leq t \leq T.
\end{equation}
\end{thm}

\begin{cor}[Convergence of interface]\label{total-2-chemo}
There exists $C>0$ such that
\begin{equation}\label{thi-3-chemo}
d_\mathcal H (\Gamma ^\ep _t,\Gamma _t)\leq C\ep \quad\ \text
{for}\; \  0\leq t \leq T,
\end{equation}
where
$$
d_\mathcal H (A,B):=\max \{ \sup_{a \in A} d(a,B),\,\sup_{b\in B}
d(b,A)\}
$$
denotes the Hausdorff distance between two compact sets $A$ and
$B$.  Consequently, $\Gamma^\ep_t \to \Gamma_t$ as $\ep\to 0$,
uniformly in $0\leq t\leq T$, in the sense of the Hausdorff
distance.
\end{cor}

As far as we know, the best thickness estimate in the literature
was of order $\ep |\ln \ep|$ (see \cite{C1}, \cite{C2}). We refer
to a forthcoming article \cite{KNHM}, respectively \cite{AHM}, in
which an order $\ep$ estimate is established for a Lotka-Volterra
competition-diffusion system, respectively for the FitzHugh-Nagumo
system.

\vskip 8 pt The organization of this paper is as follows. Section
\ref{preliminaries-chemo} is devoted to preliminaries: we recall
the method of asymptotic expansions to derive the equation of the
interface motion; we also recall a relaxed comparison principle
used in \cite{BHLM3}. In Section \ref{generation-chemo}, we prove
a generation of interface property. The corresponding sub- and
super-solutions are constructed by modifying the solution of the
ordinary differential equation $u_t=\ep^{-2}f(u)$, obtained by
neglecting diffusion and chemotaxis. In Section
\ref{motion-chemo}, in order to study the motion of interface, we
construct a pair of sub- and super-solutions that rely on a
related one-dimensional stationary problem. Finally, in Section
\ref{s:proof-chemo}, by fitting the two pairs of sub- and
super-solutions into each other, we prove Theorem
\ref{width-chemo}, Theorem \ref{error-chemo} and theirs
corollaries.

\section{Some preliminaries}\label{preliminaries-chemo}

\subsection{Formal derivation}\label{formal-chemo}

A formal derivation of the equation of interface motion was given
in \cite{BHLM2}. Nevertheless we briefly present it in a slightly
different way: we use arguments similar to those in \cite{NMHS}
where the first two terms of the asymptotic expansion determine
the interface equation. The observations we make here will help
the rigorous analysis in later sections, in particular for the
construction of sub- and super-solutions for the study of the
motion of interface in Section \ref{motion-chemo}.

Let $(u^\ep,v^\ep)$ be the solution of Problem $\Pe$. We recall
that $\Gamma _t ^\ep:=\{x \in \om, u^\ep(x,t)=1/2\}$ is the
interface at time $t$ and call $\Gamma ^\ep:=\bigcup _{t \geq 0}
(\Gamma _t ^\ep \times \{t\})$ the interface. Let $\Gamma=\bigcup
_{0\leq t \leq T}(\Gamma_t\times\{t\})$ be the solution of the
limit geometric motion problem and let $\widetilde d$ be the
signed distance function to $\Gamma$ defined by:
\begin{equation}\label{eq:dist-chemo}
\widetilde d (x,t)=
\begin{cases}
&\mbox{dist}(x,\gam _t) \quad \text{for }x\in\om _t^{(0)} \\
-&\mbox{dist}(x,\gam _t) \quad \text{for } x\in\om _t^{(1)} ,
\end{cases}
\end{equation}
where $\mbox{dist}(x,\Gamma _t)$ is the distance from $x$ to the
hypersurface  $\Gamma _t$ in $\om$. We remark that $\widetilde
d=0$ on $\Gamma$ and that $|\nabla \widetilde d|=1$ in a
neighborhood of $\Gamma.$ We then define
\begin{equation*}
Q^{(1)}_T = \bigcup_{\,0<t\leq T}(\Omega^{(1)}_t
\times\{t\}),\qquad Q^{(0)}_T = \bigcup_{\,0<t\leq
T}(\Omega^{(0)}_t \times\{t\}).
\end{equation*}
We assume that the solution $\ue$ has the expansions
\begin{equation} \label{outer-chemo}
\ue(x,t)= \{0\;\text{or}\; 1\}  + \ep u_1 (x,t) + \cdots
\end{equation}
away from the interface $\Gamma$ (the outer expansion) and
\begin{equation}\label{inner-chemo}
\ue(x,t)=\U(x,t,\frac{\widetilde d(x,t)}{\ep})+\ep
\V(x,t,\frac{\widetilde d(x,t)}{\ep})+\cdots
\end{equation}
near $\Gamma$ (the inner expansion). Here, the functions
$U_k(x,t,z)$, $k=0,1,\cdots$, are defined for $x\in \overline
\Omega$, $t\geq 0$, $z\in \R$. The stretched space variable
$\xi:=\widetilde d(x,t)/\ep$ gives exactly the right spatial
scaling to describe the rapid transition between the regions
$\{\ue \approx 1\}$ and $\{\ue \approx 0\}$. We use the
normalization conditions
$$
U_0(x,t,0)=1/2, \quad U_k(x,t,0)=0,
$$
for all $k \geq 1$.  The matching conditions between the outer and
the inner expansion are given by
\begin{equation}\label{match-chemo}
\begin{array}{ll}
U_0(x,t,+\infty)=0, \quad &U_k(x,t,+\infty)= 0, \vspace{3pt}\\
U_0(x,t,-\infty)=1, \quad &U_k(x,t,-\infty)= 0,
\end{array}
\end{equation}
for all $k \geq 1$. We also assume that the solution $\ve$ has the
expansion
\begin{equation}\label{v-expansion-chemo}
\ve(x,t)= v_0(x,t)  + \ep v_1 (x,t)+ \cdots
\end{equation}
in $\om \times (0,T)$.

We now substitute the inner expansion \eqref{inner-chemo} and the
expansion \eqref{v-expansion-chemo} into the parabolic equation of
$\Pe$ and collect the $\ep ^{-2}$ terms. We omit the calculations
and, using $|\nabla\tilde{d}|=1$ near $\Gamma_t$, the
normalization and matching conditions, we deduce that
$U_0(x,t,z)=U_0(z)$ is the unique solution of the stationary
problem
\begin{equation}\label{eq-phi-chemo}
\left\{\begin{array}{ll}
{U_0} '' +f(U_0)=0 \vspace{3pt}\\
U_0(-\infty)= 1,\quad U_0(0)=1/2,\quad U_0(+\infty)= 0.
\end{array} \right.
\end{equation}
This solution represents the first approximation of the profile of
a transition layer around the interface observed in the stretched
coordinates. Recalling that the nonlinearity is given by
$f(u)=u(1-u)(u-1/2)$, we have
\begin{equation}\label{formulepouruzero-chemo}
U_0(z)=\frac 12\big(1-\tanh \frac {z}{2\sqrt
2}\big)=\frac{e^{-z/\sqrt 2}}{1+e^{-z/\sqrt 2}}\ .
\end{equation}
We claim that $U_0$ has the following properties.

\begin{lem}\label{est-phi-chemo}
There exist positive constants $C$ and $\lambda$ such that the
following estimates hold.
$$
\begin{array}{ll}
0 <\, \, \, \,U_0(z)&\leq Ce^{-\lambda|z|} \quad \text{ for } z\geq 0,\vspace{3pt}\\
0 <1-U_0(z)&\leq Ce^{-\lambda|z|} \quad \text{ for } z\leq 0.\\
\end{array}
$$
In addition, $U_0$ is a strictly decreasing function and
$$
|{U_0}'(z)|+|{U_0}''(z)|\leq Ce^{-\lambda|z|} \quad \text{ for }
z\in \mathbb{R}.
$$
\end{lem}

\vskip 8 pt Next we collect the $\ep ^{-1}$ terms. Since $U_0$
depends only on the variable $z$, we have $\n U_{0z}=0$ which,
combined with the fact that $|\nabla\tilde{d}|=1$ near $\Gamma_t$,
yields
\begin{equation}\label{eqU1-chemo}
\VVV + f'(\U)\V ={U_0}'(\widetilde d_{t}-\Delta \widetilde d+\n
\widetilde d \cdot \n \chi(v^0))- g(\U),
\end{equation}
a linearized problem corresponding to \eqref{eq-phi-chemo}. The
solvability condition for the above equation, which can be seen as
a variant of the Fredholm alternative, plays the key role for
deriving the equation of interface motion. It is is given by
$$
\int_\BR \Big[ {{U_0}'}^2(z)(\widetilde d_{t}-\Delta \widetilde
d+\n \widetilde d \cdot \n \chi(v^0))(x,t) - g(U_0(z)){U_0}'(z)
\Big ]dz=0,
$$
for all $(x,t) \in Q_T$. By the definition of $g$ in
\eqref{nonlinearity-chemo}, we compute
$$
\int _\R g(U_0(z)){U_0}'(z)dz=-\int _0 ^1 g(u)du=-\alpha /6,
$$
whereas the equality \eqref{formulepouruzero-chemo} yields
$$
\int _\R {U_0}'^2(z)dz=\frac{1}{\sqrt 2} \int _0 ^{+\infty}
\frac{u}{1+u^4}du=1/6\sqrt 2.
$$
Combining the above expressions, we obtain
\begin{equation}\label{eq-d-chemo}
\Big(\widetilde d_{t}-\Delta \widetilde d+\n \widetilde d \cdot \n
\chi(v^0)\Big)(x,t)=-\sqrt 2 \alpha.
\end{equation}
Since $\n \widetilde d$ $(=\nabla_x\, \widetilde d(x,t))$
coincides with the outward normal unit vector to the hypersurface
$\Gamma _t$, we have $\widetilde d_{t}(x,t)=-V_n$, where $V_n$ is
the normal velocity of the interface $\Gamma _t$. It is also known
that the mean curvature $\kappa$ of the interface is equal to
$\Delta \widetilde d/(N-1)$. Thus the above equation reads as
\begin{equation}\label{eq:2-6-chemo}
V_{n}=-(N-1)\kappa + \di{\frac{\partial {\chi(v^0)}}{\partial n}}
+ \sqrt 2 \alpha \quad \ \textrm{on}\ \; \gam _t,
\end{equation}
that is the equation of interface motion in $(P^0)$. Summarizing,
under the assumption that the solution $\ue$ of Problem $\Pe$
satisfies
\begin{equation*}
\ue\to
\begin{cases}
1 &\quad  \textrm{ in } Q_T^{(1)} \\
0 &\quad  \textrm{ in } Q_T^{(0)}
\end{cases}\qquad\hbox{as}\ \ \ep\to 0,
\end{equation*}
we have formally showed that the boundary $\Gamma _t$ between
$\Omega_t^{(0)}$ and $\Omega_t^{(1)}$ moves according to the law
\eqref{eq:2-6-chemo}.

One can note that, using the equality
\eqref{formulepouruzero-chemo}, we clearly have $\sqrt 2 \alpha
{U_0}'+g(U_0)\equiv 0$ so that, substituting \eqref{eq-d-chemo}
into \eqref{eqU1-chemo} yields $U_1\equiv 0$.

\subsection{A comparison principle}
The definition of sub- and super-solutions is the one proposed in
\cite{BHLM3}.

\begin{defi}\label{def-sub-sup-chemo}
Let $(u_\ep ^-,u_ \ep ^+)$ be two smooth functions with $u_\ep
^-\leq u_ \ep ^+$ on $Q_T$ and
$$
\frac{\partial u_\ep ^-}{\partial \nu} \leq 0 \leq \frac{\partial
u_\ep ^+}{\partial \nu} \quad\ \text{on} \; \ \partial \om \times
(0,T).
$$
By definition, $(u_\ep ^-,u_ \ep ^+)$ is a pair of sub- and
super-solutions if, for any $v^\ep$ which satisfies
\begin{equation}\label{v-coincee-def-chemo}
\begin{cases}
u_\ep ^- \leq -\Delta v^\ep +\gamma v^\ep \leq u_ \ep ^+ \quad &\text{ on } Q_T,\vsp \\
\di{\frac{\partial v^\ep}{\partial \nu}}=0 \quad &\text{ on }
\partial \om\times (0,T),
\end{cases}
\end{equation}
we have
$$
L_{v^\ep}[u_\ep ^-] \leq 0 \leq L_{v^\ep}[u_\ep ^+],
$$
where the operator $L_{v^\ep}$ is defined by
$$
L_{v^\ep}[\phi]=\phi _t-\Delta \phi +\n \cdot (\phi\n
\chi(v^\ep))-\di{\frac{1}{\ep ^2}}\fe(\phi).
$$
\end{defi}

\vskip 8 pt As proved in \cite{BHLM3}, the following comparison
principle holds.

\begin{prop}\label{comparison-chemo}
Let a pair of sub- and super-solutions be given. Assume that, for
all $x \in \om$,
$$
u_\ep ^-(x,0) \leq u_0(x) \leq u_\ep ^+(x,0).
$$
Then, if we denote by $(\ue,\ve)$ the solution of Problem $\Pe$,
the function $\ue$ satisfies, for all $(x,t)\in Q_T$,
$$
u_\ep ^-(x,t) \leq \ue (x,t) \leq u_\ep ^+(x,t).
$$
\end{prop}

\section{Generation of interface}\label{generation-chemo}

In this section we study the rapid formation of internal layers in
a neighborhood of $\Gamma_0=\{x\in \om,u_0(x)=1/2\}$ within a very
short time interval of order $\ep^2 |\ln\ep|$. In the sequel, we
shall always assume that $0<\eta < 1/4$. The main result of this
section is the following.

\begin{thm}\label{th-gen-chemo}
Let $\eta$ be arbitrary and define $\mu$ as the derivative of
$f(u)$ at the unstable equilibrium $u=1/2$, that is
\begin{equation}\label{def-mu-chemo}
\mu=f'(1/2)=1/4.
\end{equation}
Then there exist positive constants $\ep_0$ and $M_0$ such that,
for all $\,\ep \in (0,\ep _0)$,
\begin{itemize}
\item for all $x\in\om$,
\begin{equation}\label{part1-chemo}
-\eta \leq u^\ep(x,\mu ^{-1} \ep ^2 | \ln \ep |) \leq 1+\eta,
\end{equation}
\item for all $x\in\om$ such that $|u_0(x)-\frac 12|\geq M_0 \ep$,
we have that
\begin{align}
&\text{if}\;~~u_0(x)\geq 1/2+M_0\ep\;~~\text{then}\;~~u^\ep(x,\mu
^{-1} \ep ^2 | \ln \ep |)
\geq 1-\eta,\label{part2-chemo}\\
&\text{if}\;~~u_0(x)\leq 1/2-M_0\ep\;~~\text{then}\;~~u^\ep(x,\mu
^{-1} \ep ^2 | \ln \ep |)\leq \eta \label{part3-chemo}.
\end{align}
\end{itemize}
\end{thm}

\vskip 8 pt The above theorem will be proved by constructing a
suitable pair of sub and super-solutions.

\subsection{The perturbed bistable ordinary differential equation}

\vskip 8 pt We first consider a slightly perturbed nonlinearity:
$$
\f(u)=f(u)+\delta ,
$$
where $\de$ is any constant. For $|\delta|$ small enough, this
function is still cubic and bistable; more precisely, we claim
that it has the following properties.

\begin{lem}Let $\de _0>0$ be small enough. Then, for all $\de \in
(-\de _0,\de _0)$,
\begin{itemize}
\item $\f$ has exactly three zeros, namely $\hm<\h<\hp$. More
precisely,
\begin{equation}\label{f-delta-chemo}
\f (u)=(u-\am (\de))(\ap (\de)-u)(u-a(\de)),
\end{equation}
and there exists a positive constant $C$ such that
\begin{equation}\label{h-chemo}
|\hm |+|\h -1/2|+|\hp-1|\leq C|\de|.
\end{equation}
\item We have that
\begin{equation}\label{signeF-chemo}
\begin{array}{ll}\f \quad \text { is strictly
positive in }\quad  (-\infty,\hm)\cup(\h,\hp),\vspace{3pt}\\
\f \quad \text{ is strictly negative in } \quad
(\hm,\h)\cup(\hp,+\infty).
\end{array}
\end{equation}
\item Set
$$
\m:= \f '(\h)=f'(\h),
$$
then there exists a positive constant, which we denote again by
$C$, such that
\begin{equation}\label{mu-chemo}
| \m -\mu| \leq C|\de|.
\end{equation}
\end{itemize}
\end{lem}

\vskip 8 pt In order to construct a pair of sub and
super-solutions for Problem $\Pe$ we define $Y(\tau,\xi;\de)$ as
the solution of the ordinary differential equation
\begin{equation}\label{ode-chemo}
\left\{\begin{array}{ll} Y_\tau (\tau,\xi;\de)&=\f
(Y(\tau,\xi;\de)) \quad\ \text {for} \; \ \tau >0,\vspace{3pt}\\
Y(0,\xi;\de)&=\xi ,
\end{array}\right.
\end{equation}
for $\de \in (-\de _0,\de _0)$ and $\xi\in (-2C_0,2C_0)$, where
$C_0$ has been chosen in \eqref{int1-chemo}. We present below
basic properties of $Y$.

\begin{lem}\label{Y1-chemo}
We have $Y_\xi >0$, for all $\xi \in (-2C_0,2C_0)\setminus
\{\am(\de),a(\de),\ap(\de)\}$, all $\de\in(-\de _0,\de _0)$ and
all $\tau
> 0$. Furthermore,
$$
\Yxi=\frac{\fde(\Y)}{\fde(\xi)}.
$$
\end{lem}

{\noindent \bf Proof.} We differentiate \eqref{ode-chemo} with
respect to $\xi$ to obtain
\begin{equation*}
\left\{\begin{array}{ll} Y_{\xi\tau}=Y_\xi f'(Y),\vspace{3pt}\\
Y_\xi(0,\xi;\de)=1,
\end{array}\right.
\end{equation*}
which is integrated as follows:
\begin{equation}\label{Y2-chemo}
\Yxi=\exp \Big[\int_0^\tau f'(Y(s,\xi;\de))ds\Big]>0.
\end{equation}
Then differentiating \eqref{ode-chemo} with respect to $\tau$, we
obtain
\begin{equation*}
\left\{\begin{array}{ll} Y_{\tau\tau}=Y_\tau f'(Y),\vspace{3pt}\\
Y_\tau(0,\xi;\de)=\fde (\xi),
\end{array}\right.
\end{equation*}
which in turn implies
$$
Y_\tau(\tau,\xi;\de)=\fde(\xi) \exp \Big[\int_0^\tau
f'(Y(s,\xi;\de))ds\Big],
$$
which enables to conclude. \qed

\vskip 8pt We define a function $\A$ by
\begin{equation}\label{A-chemo}
\A=\frac{f'(\Y)-f'(\xi)}{\fde(\xi)}.
\end{equation}

\begin{lem}\label{Y5-chemo}
We have, for all $\xi \in (-2C_0,2C_0)\setminus
\{\am(\de),a(\de),\ap(\de)\}$, all $\de \in(-\de _0,\de _0)$ and
all $\tau
> 0$,
$$
\A=\int_0^\tau f''(Y(s,\xi;\de))Y_\xi(s,\xi;\de)ds.
$$
\end{lem}

{\noindent \bf Proof.} We differentiate the equality of Lemma
\ref{Y1-chemo} with respect to $\xi$ to obtain
\begin{equation}\label{Y4-chemo}
\Yxixi=\A \Yxi.
\end{equation}
Then differentiating \eqref{Y2-chemo} with respect to $\xi$ yields
$$
Y_{\xi\xi}=Y_\xi \int_0^\tau f''(Y(s,\xi;\de))Y_\xi(s,\xi;\de)ds.
$$
These two last results complete the proof of Lemma \ref{Y5-chemo}.
\qed

\vskip 8pt Next we prove estimates on the growth of $Y$, $A$ and
theirs derivatives. We first consider the case where the initial
value $\xi$ is far from the stable equilibria, more precisely when
it lies between $\eta$ and $1-\eta$.

\begin{lem}\label{est-derY-A-milieu-chemo}
Let $\eta$ be arbitrary. Then there exist positive constants
$\de_0=\de_0(\eta)$, $\widetilde C _1=\widetilde C _1(\eta)$,
$\tilde C _2=\tilde C _2(\eta)$ and $C_3=C_3(\eta)$ such that, for
all $\de \in (-\de _0,\de_0)$, for all $\tau>0$,
\begin{itemize}
\item if $\xi \in (\h,1-\eta)$ then, for every $\tau >0$ such that
$Y(\tau,\xi;\de)$ remains in the interval $(\h,1-\eta)$, we have
\begin{equation}\label{est-Y3-chemo}
\tilde C _1e^{\m\tau}\leq Y_\xi(\tau,\xi;\de) \leq \tilde C _2
e^{\m\tau},
\end{equation}
and
\begin{equation}\label{est-A-milieu-chemo}
|A(\tau,\xi;\de)|\leq C_3(e^{\m\tau}-1);
\end{equation}
\item if $\xi\in (\eta,\h)$ then, for every $\tau >0$ such that
$Y(\tau,\xi;\de)$ remains in the interval $(\eta,\h)$,
\eqref{est-Y3-chemo} and \eqref{est-A-milieu-chemo} hold as well.
\end{itemize}
\end{lem}

{\noindent \bf Proof.} We take $\xi \in (a(\de),1-\eta)$ and
suppose that for $s \in (0,\tau)$, $Y(s,\xi;\de)$ remains in the
interval $(a(\de),1-\eta)$. Integrating the equality
$$
\frac{Y_\tau(s,\xi;\de)}{\fde (Y(s,\xi;\de))}=1
$$
from $0$ to $\tau$ and using the change of variable
$q=Y(s,\xi;\de)$ leads to
\begin{equation}\label{g-tau-chemo}
\int _\xi ^{Y(\tau,\xi;\de)} \frac{dq}{\fde(q)}=\tau.
\end{equation}
Moreover, the equality in Lemma \ref{Y1-chemo} enables to write
\begin{equation}\label{naka-chemo}
\begin{array}{lll}
\ln Y_\xi (\tau,\xi;\de)&=
\di{\int _ \xi ^{Y(\tau,\xi;\de)}} \frac{f'(q)}{\fde(q)}dq \vsp \\
&=\di{\int _ \xi ^{Y(\tau,\xi;\de)}}\big
[\frac{f'(a(\de))}{\fde(q)}+\frac{f'(q)-f'(a(\de))}{\fde(q)}\big ]dq \vsp \\
&=\mu (\de) \tau+\di{\int _ \xi ^{Y(\tau,\xi;\de)}}\hde(q)dq,
\end{array}
\end{equation}
where
$$
h_\de (q)=\frac{f'(q)-f'(a(\de))}{f_\de(q)}.
$$
In view of \eqref{mu-chemo}, respectively \eqref{h-chemo}, we can
choose $\de _0=\de _0 (\eta) >0$ small enough so that, for all
$\de \in [-\de _0,\de _0]$, we have $\mu(\de)\geq \mu /2>0$,
respectively $(a(\de),1-\eta]\subset(a(\de),\ap (\de))$. Since
\[
h_\de(q)\to \frac{\f''(a(\de))}{\f'(a(\de))}
=\frac{f''(a(\de))}{f'(a(\de))} \quad\ \hbox{as} \ \ q \to a(\de),
\]
we see that the function $(q,\de)\mapsto h_\de(q)$ is continuous
in the compact region $\{\,|\de|\leq \de_0,\;a(\de)\leq q\leq
1-\eta\,\}$. It follows that $|h_\de(q)|$ is bounded by a constant
$H=H(\eta)$ as $(q,\de)$ varies in this region. Since
$|Y(\tau,\xi;\de)-\xi|$ takes its values in the interval
$[0,1-\eta-a(\de)]\subset[0,1]$, it follows from
\eqref{naka-chemo} that
$$
\mu (\de) \tau -H \leq \ln Y_\xi(\tau,\xi;\de) \leq \mu (\de)
\tau+H,
$$
which, in turn, proves \eqref{est-Y3-chemo}. Next Lemma
\ref{Y5-chemo} and \eqref{est-Y3-chemo} yield
$$\begin{array}{ll}
|A(\tau,\xi;\de)| &\leq \Vert f'' \Vert _{L^\infty (0,1)}
\di{\int_0^\tau} \tilde C _2 e^{\mu (\de) s}ds \vsp \\
&\leq \di{\frac{\Vert f'' \Vert _{L^\infty (0,1)}\tilde C
_2}{\mu(\de)}}(e^{\mu (\de) \tau}-1)\vsp\\
&\leq \di{\frac{2}{\mu}}\Vert f'' \Vert _{L^\infty (0,1)}\tilde C
_2(e^{\mu (\de) \tau}-1),
\end{array}
$$
which completes the proof of \eqref{est-A-milieu-chemo}. The case
where $\xi$ and $Y(\tau,\xi;\de)$ are in $(\eta,a(\de))$ is
similar and omitted. \qed

\begin{cor}\label{est-Y-milieu-chemo}
Let $\eta$ be arbitrary. Then there exist positive constants
$\de_0=\de_0(\eta)$, $C_1=C_1(\eta)$ and $C_2=C_2(\eta)$ such
that, for all $\de \in (-\de _0,\de _0)$, for all $\tau>0$,
\begin{itemize}
\item if $\xi\in (\h,1-\eta)$ then, for every $\tau >0$ such that
$Y(\tau,\xi;\de)$ remains in the interval $(\h,1-\eta)$, we have
\begin{equation}\label{est-Y-1-chemo}
C_1e^{\m \tau}(\xi-\h)\leq Y(\tau,\xi;\de)-\h \leq C_2e^{\m
\tau}(\xi-\h),
\end{equation}
\item if $\xi \in (\eta,\h)$ then, for every $\tau >0$ such that
$Y(\tau,\xi;\de)$ remains in the interval $(\eta,\h)$, we have
\begin{equation}\label{est-Y-2-chemo}
C_2e^{\m \tau}(\xi-\h)\leq Y(\tau,\xi;\de)-\h \leq C_1e^{\m
\tau}(\xi-\h).
\end{equation}
\end{itemize}
\end{cor}

{\noindent \bf Proof.} In view of \eqref{mu-chemo}, respectively
\eqref{h-chemo}, we can choose $\de _0=\de _0 (\eta) >0$ small
enough so that, for all $\de \in [-\de _0,\de _0]$, we have
$\mu(\de)\geq \mu /2>0$, respectively
$(a(\de),1-\eta]\subset(a(\de),\ap (\de))$. Since
\[
\frac{\f(q)}{q-a(\de)}\to \mu(\de) \quad\ \hbox{as} \ \ q \to
a(\de),
\]
it follows that $(q,\de)\mapsto\f(q)/(q-a(\de))$ is a strictly
positive and continuous function in the compact region
$\{\,|\de|\leq \de_0,\;a(\de)\leq q\leq 1-\eta\,\}$, which insures
the existence of constants $B_1=B_1(\eta)>0$ and $B_2=B_2(\eta)>0$
such that, for all $q\in(\h,1-\eta)$, all $\de\in(-\de _0,\de
_0)$,
\begin{equation}\label{ineg-fq-chemo}
B_1(q-a(\de))\leq f_ \de (q) \leq B_2(q-a(\de)).
\end{equation}
We write the inequalities \eqref{ineg-fq-chemo} for $q=\Y
\in(\h,1-\eta)$ and then for $q=\xi \in (\h,1-\eta)$, which,
together with Lemma \ref{Y1-chemo}, implies that
$$
\frac{B_1}{B_2}(\Y-\h)\leq (\xi-\h)\Yxi \leq
\frac{B_2}{B_1}(\Y-\h).
$$
In view of \eqref{est-Y3-chemo}, this completes the proof of
inequalities \eqref{est-Y-1-chemo}. The proof of
\eqref{est-Y-2-chemo} is similar and omitted.\qed

\vskip 8pt We now present estimates in the case that the initial
value $\xi$ is smaller than $\eta$ or larger than $1-\eta$.

\begin{lem}\label{est-bords-chemo}
Let $\eta$ and $M>0$ be arbitrary. Then there exist positive
constants $\de _0=\de _0(\eta,M)$ and $C_4=C_4(M)$ such that, for
all $\de \in (-\de _0,\de _0)$,
\begin{itemize}
\item if $\xi \in [1-\eta,1+M]$, then, for all $\tau > 0$,
$Y(\tau,\xi;\de)$ remains in the interval $[1-\eta,1+M]$ and
\begin{equation}\label{est-A-bords-chemo}
|A(\tau,\xi;\de)|\leq C_4\tau \quad\hbox{for}\ \ \tau>0 \,;
\end{equation}
\item if $\xi \in [-M,\eta]$, then, for all $\tau >0$,
$Y(\tau,\xi;\de)$ remains in the interval $[-M,\eta]$ and
\eqref{est-A-bords-chemo} holds as well.
\end{itemize}
\end{lem}

{\noindent \bf Proof.} Since the two statements can be treated in
the same way, we will only prove the former. The fact that $\Y$,
the solution of the ordinary differential equation
\eqref{ode-chemo}, remains in the interval $[1-\eta,1+M]$ directly
follows from the bistable properties of $\f$, or, more precisely,
from the sign conditions $\f(1-\eta)>0$, $\f(1+M)<0$ valid if $\de
_0=\de _0(\eta,M)$ is small enough.

To prove \eqref{est-A-bords-chemo}, suppose first that $\x\in
[\hp,1+M]$. By the above arguments, $\Y$ remains in this interval.
Moreover $f'$ is negative in this interval. Hence, it follows from
\eqref{Y2-chemo} that $\Yxi \leq 1$. We then use Lemma
\ref{Y5-chemo} to deduce that
$$
|\A| \leq \Vert f'' \Vert _{L^\infty (-M,1+M)}\tau=:C_4\tau.
$$
The case $\xi \in [1-\eta,\hp]$ being similar, this completes the
proof of the lemma. \qed

\vskip 8pt Now we choose the constant $M$ in the above lemma
sufficiently large so that $[-2C_0,2C_0]\subset [-M,1+M]$, and fix
$M$ hereafter. Therefore the constant $C_4$ is fixed as well.
Using the fact that $\tau\mapsto \tau(e^{\mu(\de) \tau}-1)^{-1}$
is uniformly bounded for $\de \in (-\de _0, \de _0)$, with $\de
_0$ small enough (see \eqref{mu-chemo}), and for $\tau>0$, one can
easily deduce from (\ref{est-A-milieu-chemo}) and
(\ref{est-A-bords-chemo}) the following general estimate.

\begin{lem}\label{EST-A-chemo}
Let $\eta$ be arbitrary and let $C_0$ be the constant defined in
\eqref{int1-chemo}. Then there exist positive constants
$\de_0=\de_0(\eta)$, $C_5=C_5(\eta)$ such that, for all $\de \in
(-\de _0,\de _0)$, all $\xi \in (-2C_0,2C_0)$ and all $\tau>0$,
$$
|A(\tau,\xi;\de)|\leq C_5(e^{\m\tau}-1).
$$
\end{lem}

\subsection{Construction of sub and super-solutions}
We now use $Y$ to construct a pair of sub- and super-solutions for
the proof of the generation of interface theorem. We set
\begin{equation}\label{w+--general-chemo}
w_\ep^\pm(x,t)=Y\Big(\frac{t}{\ep^2},u_0(x)\pm\ep^2r(\pm \ep
G,\frac{t}{\ep^2});\pm \ep G\Big),
\end{equation}
where the constant $G$ is defined by
$$
G=\sup_{u\in [-2C_0,2C_0]} |g(u)|,
$$
and the function $r(\de,\tau)$ is given by
$$
r(\de,\tau)=C_6(e^{\m\tau}-1).
$$
For simplicity, we make the following additional assumption:
\begin{equation}\label{int2-chemo}
\frac{\partial u_0}{\partial \nu} =0 \quad\ \textrm{on}\ \;\dom.
\end{equation}
In the general case where \eqref{int2-chemo} does not necessary
hold, we have to slightly modify $w _\ep ^\pm$ near the boundary
$\partial \om$. This will be discussed in the next remark.

\begin{lem}\label{w-chemo}
There exist positive constants  $\ep_0$ and $C_6$ such that for
all $\,\ep \in (0,\ep _0)$, the functions $w_\ep^-$ and $w_\ep^+$
are respectively sub- and super-solutions for Problem $(P^\ep)$,
in the domain
$$
 \big{\{}(x,t)\in \Q, \; x\in \om,\; 0\leq t \leq
\mu ^{-1} \ep^2|\ln \ep|\big{\}}.
$$
\end{lem}

{\noindent \bf Proof.} First, \eqref{int2-chemo} implies the
homogeneous Neumann boundary condition
$$
\frac{\partial w_\ep^\pm}{\partial \nu}=0 \quad\ \textrm{on}\
\;\partial \Omega \times (0,+\infty).
$$
Let $v^\ep$ be such that
\begin{equation}\label{v-coincee1-chemo}
\begin{cases}
w_\ep ^- \leq -\Delta v^\ep +\gamma v^\ep \leq w_ \ep ^+ \vsp \\
\di{\frac{\partial v^\ep}{\partial \nu}}=0.
\end{cases}
\end{equation}
According to Definition \ref{def-sub-sup-chemo}, what we have to
show is
$$
L_{v^\ep}[w_\ep^+]:=(w_\ep^+) _t-\Delta w_\ep^+ +\n \cdot
(w_\ep^+\n \chi(v^\ep))-\frac{1}{\ep ^2}\fe(w_\ep^+)\geq 0.
$$
Let $C_6$ be a positive constant which does not depend on $\ep$.
If $\ep _0$ is sufficiently small, we note that $\pm \ep G \in
(-\de _0,\de _0)$ and that, in the range $0 \leq t \leq \mu ^{-1}
\ep ^2|\ln \ep|$,
$$
|\ep^2 C_6(e^{\mu(\pm \ep G)t/\ep^2}-1)| \leq \ep
^2C_6(\ep^{-\mu(\pm\ep G)/\mu}-1) \leq C_0,
$$
which implies that
$$
u_0(x)\pm\ep^2r(\pm \ep G,t/\ep^2) \in (-2C_0,2C_0).
$$
These observations allow us to use the results of the previous
subsection with $\tau =t/\ep^2$, $\xi=u_0(x)+\ep^2r(\ep
G,t/\ep^2)$ and $\de=\ep  G$. In particular, setting $F_1:=\Vert
f' \Vert _{L ^\infty (-2C_0,2C_0)}$, this implies, using
\eqref{Y2-chemo}, that
$$
e^{-F_1 T} \leq Y_\xi \leq e^{F_1 T}.
$$
Straightforward computations yield
\begin{eqnarray*}
L_{v^\ep}[w_\ep^+]&=&\edeux Y_\tau + C_6\mu(\ep G) \emuttchemo
Y_\xi -\Delta u_0 Y _ \xi-|\n u_0|^2 Y_{\xi \xi}\\
& &{}+Y_\xi \n u_0 \cdot \n \chi(v^\ep)+Y \Delta
\chi(v^\ep)-\edeux f(Y)-\frac 1 \ep g(Y),
\end{eqnarray*}
and then, in view of the ordinary differential equation
\eqref{ode-chemo}, $\ep G$ playing the role of $\de$,
\begin{multline}
L_{v^\ep}[w_\ep^+]=\di{\frac {1}{\ep}}\big[
G-g(Y)\big]+Y_\xi\Big[C_6 \mu(\ep
 G)\emuttchemo-\Delta u_0-\displaystyle{\frac{Y_{\xi\xi}}{Y_\xi}}|\nabla u_0|^2\vsp\\
 +\n u_0 \cdot \n \chi (v^\ep)+\di{\frac{Y}{Y_\xi}}\Delta
\chi(v^\ep)\Big].
\end{multline}
By the definition of $G$ the first term above is positive. Now,
using the choice of $C_0$ in \eqref{int1-chemo}, the fact that
$Y_{\xi \xi}/Y_ \xi=A$ and Lemma \ref{EST-A-chemo}, we obtain, for
a $C_5$ independent of $\ep$,
\begin{eqnarray*}
L_{v^\ep}[w_\ep^+]\geq & Y_\xi\Big[C_6 \mu(\ep
G)\emuttchemo-C_0-C_5(\emuttchemo -1)C_0 ^2 \\
& {}{}{}-C_0 |\n \chi (v^\ep)|-2C_0 e^{F_1 T} |\Delta \chi
(v^\ep)| \Big ].
\end{eqnarray*}
Moreover, the inequalities in \eqref{v-coincee1-chemo} can be
written as $-\Delta v^\ep + \gamma v^\ep=h^\ep$, with $-2C_0 \leq
h^\ep \leq 2C_0$, so that the standard theory of elliptic
equations gives a uniform bound $M$ for $|v^\ep|$, $|\n v^\ep|$
and $|\Delta v^\ep|$. Hence, using the smoothness of $\chi$, we
have a uniform bound $M'$ for $|\n \chi(v^\ep)|$ and $|\Delta
\chi(v^\ep)|$. It follows that
$$
L_{v^\ep}[w_\ep^+]\geq Y_\xi \Big [(C_6 \mu(\ep G)-C_5 C_0
^2)\emuttchemo-C_0+C_5C_0 ^2-C_0 M'-2C_0 e^{F_1 T}M' \Big ].
$$
Hence, in view of \eqref{mu-chemo}, we have, for $\ep _0$ small
enough (recall that $Y_\xi
>0$),
$$
Lw_\ep^+ \geq Y_\xi\Big[(C_6\frac 1 2\mu -C_5 C_0
^2)-C_0-C_0M'-2C_0 e^{F_1 T}M' \Big ] \geq 0,
$$
for $C_6$ large enough, so that $w_\ep^+$ is a super-solution for
Problem $\Pe$. We omit the proof that $w_\ep^-$ is a sub-solution.
\qed

\vskip 8pt Now, since $w_\ep^\pm(x,0)=Y(0,u_0(x);\pm \ep
G)=u_0(x)$, the comparison principle set in Proposition
\ref{comparison-chemo} asserts that, for all $x \in \Omega$, for
all $0\leq t \leq \mu ^{-1} \ep ^2 |\ln \ep|$,
\begin{equation}\label{coincee-chemo}
w_\ep^-(x,t) \leq u^\ep(x,t) \leq w_\ep^+(x,t).
\end{equation}

\begin{rem} In the more general case where \eqref{int2-chemo} is not
valid, one can proceed in the following way: in view of
\eqref{dalltint-chemo} and \eqref{initial-data-chemo} there exist
positive constants $d_1$ and $\rho$ such that $u_0(x) \leq
1/2-\rho$ if $d(x,\partial \Omega) \leq d_1$. Let $\chi$ be a
smooth cut-off function defined on $[0,+\infty)$ such that $0\leq
\chi \leq 1$, $\chi(0)=\chi '(0)=0$ and $\chi (z)=1$ for $z\geq
d_1$. Then define
$$
u_0 ^+(x)=\chi(d(x,\partial \Omega))u_0(x)+(1-\chi(d(x,\partial
\Omega))(1/2 -\rho),
$$
$$
u_0 ^-(x)=\chi(d(x,\partial \Omega))u_0(x)+(1-\chi(d(x,\partial
\Omega))\min_{x\in\overline \Omega} u_0(x).
$$
Clearly, $u_0 ^- \leq u_0 \leq u_0^+$, and both $u_0^ \pm$ satisfy
\eqref{int2-chemo}. Now we set
\begin{equation*}\label{w+-2-chemo}
\tilde{w}_\ep^\pm(x,t)=
Y\Big(\frac{t}{\ep^2},\,u^\pm_0(x)\pm\ep^2r(\pm \ep
G,\frac{t}{\ep^2});\pm \ep G\Big).
\end{equation*}
Then the same argument as in Lemma \ref{w-chemo} shows that
$(\tilde{w}_\ep^-,\tilde{w}_\ep^+)$ is a pair of sub and
super-solutions for Problem $(P^\ep)$.  Furthermore, since
$\tilde{w}_\ep^-(x,0)=u_0^-(x)\leq u_0(x)\leq u_0^+(x)
=\tilde{w}_\ep^+(x,0)$, Proposition \ref{comparison-chemo} asserts
that, for all $x \in \Omega$, for all $0\leq t \leq \mu ^{-1} \ep
^2 |\ln \ep|$, we have $\tilde{w}_\ep^-(x,t) \leq u^\ep(x,t) \leq
\tilde{w}_\ep^+(x,t)$. \hfill $\square$
\end{rem}

\subsection{Proof of Theorem \ref{th-gen-chemo}}
In order to prove Theorem \ref{th-gen-chemo} we first present a
key estimate on the function $Y$ after a time of order $\tau\sim
|\ln \ep|.$

\begin{lem}
Let $\eta $ be arbitrary; there exist positive constants
$\ep_0=\ep _0(\eta)$ and $C_7=C_7(\eta)$ such that, for all $\ep
\in (0,\ep_0)$,
\begin{itemize}
\item for all $\xi\in (-2C_0,2C_0)$,
\begin{equation}\label{part11-chemo}
-\eta \leq Y(\mu ^{-1} | \ln \ep |,\xi;\pm \ep  G) \leq 1+\eta,
\end{equation}
\item for all $\xi\in (-2C_0,2C_0)$ such that $|\xi-\frac 12|\geq
C_7 \ep$, we have that
\begin{align}
&\text{if}\;~~\xi\geq 1/2+C_7 \ep\;~~\text{then}\;~~Y(\mu ^{-1}|
\ln \ep |,\xi;\pm \ep  G)
\geq 1-\eta,\label{part22-chemo}\vspace{3pt}\\
&\text{if}\;~~\xi\leq 1/2-C_7 \ep\;~~\text{then}\;~~Y(\mu ^{-1}|
\ln \ep |,\xi;\pm \ep  G)\leq \eta \label{part33-chemo}.
\end{align}
\end{itemize}
\end{lem}

{\noindent \bf Proof.} We first prove \eqref{part22-chemo}. In
view of \eqref{h-chemo}, we have, for $C_7$ large enough,
$1/2+C_7\ep \geq a(\pm \ep  G)+\frac 1 2 C_7 \ep$, for all $\ep
\in (0,\ep _0)$, with $\ep _0$ small enough. Hence for $\xi \geq
1/2+C_7 \ep$, as long as $Y(\tau,\xi;\pm \ep G)$ has not reached
$1-\eta$, we can use \eqref{est-Y-1-chemo} to deduce that
$$
\begin{array}{llll}
Y(\tau,\xi;\pm \ep G) & \geq a(\pm \ep G) +C_1 e^{\mu(\pm \ep G)
\tau}(\xi-a(\pm \ep G))\vsp \\
& \geq a(\pm \ep G) + \frac 1 2 C_1 C_7 \ep e^{\mu(\pm \ep G) \tau}\vsp \\
& \geq \frac 12 -\ep C G+  \frac 1 2 C_1 C_7 \ep e^{\mu(\pm \ep G)
\tau}\vsp \\
& \geq 1-\eta
\end{array}
$$
provided that
\begin{equation*}
\tau \geq \tau ^\ep :=\frac{1}{\mu (\pm \ep G)}\ln
\frac{1/2-\eta+C G \ep}{ C_1 C_7 \ep /2}.
\end{equation*}
To complete the proof of \eqref{part22-chemo} we must choose $C_7$
so that $\mu ^{-1}|\ln \ep|-\tau ^\ep \geq 0$. A simple
computation shows that
$$
\begin{array}{ll}
\mu ^{-1}|\ln \ep|-\tau ^\ep =\di{\frac {\mu(\pm \ep
G)-\mu}{\mu(\pm \ep G)\mu}}|\ln \ep|&-\di{\frac {1}{\mu(\pm \ep
G)}}\ln \frac{1/2-\eta +C
G\ep} {C_1/2}\vsp\\
&+\di{\frac {1}{\mu(\pm \ep G)}}\ln C_7.
\end{array}
$$
Thanks to \eqref{mu-chemo}, as $\ep\to 0$, the first term above is
of order $\ep|\ln \ep|$ and the second one of order $1$. Hence,
for $C_7$ large enough, the quantity $\mu ^{-1}|\ln \ep|-\tau
^\ep$ is positive, for all $\ep\in (0,\ep _0)$, with $\ep _0$
small enough. The proof of \eqref{part33-chemo} is similar and
omitted.

Next we prove \eqref{part11-chemo}. First note that, by taking
$\ep _0$ small enough, the stable equilibria of $f_{\pm \ep G}$,
namely $\am (\pm \ep G)$ and $\ap(\pm \ep G)$, are in
$[-\eta,1+\eta]$. Hence, $f_{\pm \ep G}$ being a bistable
function, if we leave from a $\xi \in [-\eta,1+\eta]$ then
$Y(\tau,\xi;\pm \ep  G)$ will remain in the interval
$[-\eta,1+\eta]$. Now suppose that $1+\eta \leq \xi \leq 2C_0$
(note that this work is useless if $2C_0<1+\eta$). We check below
that $Y(\mu ^{-1} | \ln \ep|,\xi;\pm \ep  G)\leq 1+\eta$. As long
as $1+\eta \leq Y \leq 2C_0$, \eqref{ode-chemo} leads to the
inequality $Y_\tau \leq f(1+\eta)+\ep  G \leq \frac 1 2
f(1+\eta)<0$, for $\ep_0=\ep_0 (\eta)$ small enough. By
integration from $0$ to $\tau$, it follows that
$$
\begin{array}{lll}
Y(\tau,\xi;\pm \ep  G) &\leq \xi+\frac 1 2 f(1+\eta) \tau \vsp \\
& \leq 2C_0 +\frac 1 2 f(1+\eta) \tau \vsp \\
& \leq 1+ \eta,
\end{array}
$$
provided that
$$
\tau \geq \frac {2C_0-1-\eta}{- f(1+\eta)/2},
$$
and a fortiori for $\tau=\mu ^{-1} |\ln \ep|$, which completes the
proof of \eqref{part11-chemo}. \qed

\vskip 8pt We are now ready to prove Theorem \ref{th-gen-chemo}.
By setting $t=\mu ^{-1} \ep ^2|\ln \ep|$ in \eqref{coincee-chemo},
we obtain
\begin{multline}\label{gr-chemo}
Y\Big(\mu ^{-1}|\ln \ep|, u_0(x)-\ep^2 r(-\ep  G, \mu ^{-1}|\ln
\ep|);-\ep G\Big)\\
\leq u^\ep(x,\mu ^{-1} \ep^2|\ln \ep|) \leq Y\Big(\mu ^{-1}|\ln
\ep|, u_0(x)+\ep^2 r(\ep  G, \mu ^{-1}|\ln \ep|);+\ep  G \Big).
\end{multline}
In view of \eqref{mu-chemo},
\begin{equation}\label{point-chemo}
\lim _{\ep \rightarrow 0} \frac{\mu-\mu(\pm\ep  G)}{\mu}\ln \ep=0,
\end{equation}
so that, for $\ep _0$ small enough, we have
$$
\ep ^2r(\pm \ep  G,\mu ^{-1} |\ln \ep|)= C_6\ep(\ep
^{(\mu-\mu(\pm\ep  G))/\mu}-\ep) \in (\frac 1 2 C_6\ep, \frac 3 2
C_6 \ep).
$$
It follows that $u_0(x)\pm\ep^2r(\pm \ep  G,\mu ^{-1} |\ln \ep|)
\in (-2C_0,2C_0)$. Hence the result \eqref{part1-chemo} of Theorem
\ref{th-gen-chemo} is a direct consequence of \eqref{part11-chemo}
and \eqref{gr-chemo}.

Next we prove \eqref{part2-chemo}. We take $x\in \om$ such that
$u_0(x)\geq 1/2+M_0 \ep$ so that
$$
\begin{array}{ll}u_0(x)-\ep^2r(-\ep  G,
\mu ^{-1}|\ln \ep|)
&\geq 1/2+M_0\ep-\frac 3 2 C_6 \ep\vsp\\
&\geq 1/2+C_7 \ep,
\end{array}
$$
if we choose $M_0$ large enough. Using \eqref{gr-chemo} and
\eqref{part22-chemo} we obtain \eqref{part2-chemo}, which
completes the proof of Theorem \ref{th-gen-chemo}.\qed

\section{Motion of interface}\label{motion-chemo}

We have seen in Section \ref{generation-chemo} that, after a very
short time, the solution $\ue$ develops a clear transition layer.
In the present section, we show that it persists and that its law
of motion is well approximated by the interface equation in
$(P^0)$ obtained by formal asymptotic expansions in subsection
\ref{formal-chemo}.

More precisely, take the first term of the formal asymptotic
expansion \eqref{inner-chemo} as a formal expansion of the
solution:
\begin{equation}
\ue(x,t)\,\approx\,\tilde{u}^\ep(x,t):= \U\Big(\frac{\widetilde
d(x,t)}{\ep}\Big),
\end{equation}
where $U_0$ is defined in \eqref{eq-phi-chemo}. The right-hand
side is a function having a well-developed transition layer, and
its interface lies exactly on $\gam _t$. We show that this
function is a good approximation of the solution; more precisely:
\begin{quote}
{\it If $u^\ep$ becomes very close to $\tilde{u}^\ep$ at some time
moment $t=t_0$, then it stays close to $\tilde{u}^\ep$ for the
rest of time. Consequently, $\Gamma^\ep_t$ evolves roughly like
$\Gamma_t$.}
\end{quote}

To that purpose, we will construct a pair of sub- and
super-solutions $u_\ep^-$ and $u_\ep^+$ for Problem $(P^\ep)$ by
slightly modifying $\tilde{u}^\ep$. It then follows that, if the
solution $\ue$ satisfies
\[
u_\ep^-(x,t_0)\leq \ue(x,t_0)\leq  u_\ep^+(x,t_0),
\]
for some $t_0\geq 0$, then
\[
u_\ep^-(x,t)\leq \ue(x,t)\leq  u_\ep^+(x,t),
\]
for $t_0\leq t\leq T$. As a result, since both $u_\ep^+, u_\ep^-$
stay close to $\tilde{u}^\ep$, the solution $\ue$ also stays close
to $\tilde{u}^\ep$ for $t_0\leq t\leq T$.

\subsection{Construction of sub- and super-solutions}
To begin with we present mathematical tools which are essential
for the construction of sub and super-solutions.

\vskip 8 pt {\noindent \bf A modified signed distance function.}
Rather than working with the usual signed distance function
$\widetilde d$, defined in \eqref{eq:dist-chemo}, we define a
``cut-off signed distance function" $d$ as follows. Choose $d_0>0$
small enough so that $\widetilde d(\cdot,\cdot)$ is smooth in the
tubular neighborhood of $\Gamma$
$$
\{(x,t) \in \overline{\Q},\;|\widetilde{d}(x,t)|<3d_0\},
$$
and such that
\begin{equation}\label{front-chemo}
\mbox{dist}(\Gamma_t,\partial \Omega)>4d_0 \quad \textrm{ for all
} t\in [0,T].
\end{equation}
Next let $\zeta(s)$ be a smooth increasing function on $\R$ such
that
\[
\zeta(s)= \left\{\begin{array}{ll}
s &\textrm{ if }\ |s| \leq 2d_0\vspace{4pt}\\
-3d_0 &\textrm{ if } \ s \leq -3d_0\vspace{4pt}\\
3d_0 &\textrm{ if } \ s \geq 3d_0.
\end{array}
\right.
\]
We define the cut-off signed distance function $d$ by
\begin{equation}
d(x,t)=\zeta\big(\tilde{d}(x,t)\big).
\end{equation}
Note that $|\nabla d|=1$ in the region $\{(x,t) \in
\overline{\Q},\,|d(x,t)|<2d_0\}$ and that, in view of the above
definition, $\n d=0$ in a neighborhood of $\partial \Omega$. Note
also that the equation of motion interface in $(P^0)$, which is
equivalent to \eqref{eq-d-chemo}, is now written as
\begin{equation}\label{FBP-chemo}
d_t=\Delta d -\n d \cdot \n\chi(v^0)-\sqrt 2 \alpha \quad\
\textrm{on}\ \; \Gamma _t.
\end{equation}

\vskip 8 pt {\noindent \bf Construction.} We look for a pair of
sub- and super-solutions $u_\ep^{\pm}$ for Problem $\Pe$ of the
form
\begin{equation}\label{sub-chemo}
u_\ep^{\pm}(x,t)=U_0\Big(\frac{d(x,t) \mp \ep p(t)}{\ep}\Big)\pm
q(t),
\end{equation}
where $U_0$ is the solution of \eqref{eq-phi-chemo}, and where
\begin{equation}\label{defpetq-chemo}
\begin{array}{ll}
p(t)=-e^{-\beta t/\ep ^2}+e^{Lt}+ K,\vsp \\
q(t)=\sigma(\beta e^{-\beta t/\ep ^2}+\ep^2Le^{Lt}).
\end{array}
\end{equation}
Note that $q=\sigma\ep^2\,p_t$. Let us remark that the
construction \eqref{sub-chemo} is more precise than the procedure
of only taking a zeroth order term of the form $U_0$, since we
have shown in the formal derivation that the first order term
$U_1$ in \eqref{inner-chemo} vanishes. It is clear from the
definition of $u_\ep^\pm$ that
\begin{equation}\label{sub2-chemo}
\lim_{\ep\rightarrow 0} u_\ep^\pm(x,t)= \left\{
\begin{array}{ll}
1 &\textrm { for all } (x,t) \in Q_T^{(1)} \vspace{4pt}\\
0 &\textrm { for all } (x,t) \in Q_T^{(0)}.\\
\end{array}\right.
\end{equation}

The main result of this section is the following.

\begin{lem}\label{fix-chemo}
There exist positive constants $\beta,\,\sigma$ with the following
properties.  For any $K>1$, we can find positive constants $\ep_0$
and $L$ such that, for any $\ep\in(0,\ep _0)$, the functions
$u_\ep^-$ and $u_\ep^+$ are respectively sub- and super-solutions
for Problem $(P^\ep)$ in the range $x\in\ombar$, $0\leq t\leq T$.
\end{lem}

\subsection{Proof of Lemma \ref{fix-chemo}}
First, since $\n d=0$ in a neighborhood of $\partial \Omega $, we
have the homogeneous Neumann boundary condition
$$
\di{\frac{\partial u_\ep^\pm}{\partial \nu}}=0 \quad\ \textrm{on}
\ \; \partial \Omega \times [0,T].$$ Let $v^\ep$ be such that
\eqref{v-coincee-def-chemo} holds. We have to show that
$$
L_{v^\ep}[u_\ep^+]:=(u_\ep^+)_t -\Delta u_\ep^+ +\n u_\ep^+ \cdot
\n \chi(v^\ep)+u_\ep^+ \Delta \chi(v^\ep) -\edeux \fe(u_\ep^+)\geq
0,
$$
the proof of inequality $L_{v^\ep}[u_\ep^-]\leq 0$ following by
the same arguments.

\vskip 8 pt \noindent {\bf Computation of $L_{v^\ep}[u_\ep^+$].}
By straightforward computations we obtain the following terms:
$$
\begin{array}{lll}
(u_\ep^+)_t= {U_0}'(\displaystyle{\frac{d_t}{\ep}}-p_t)+q_t,\vsp \\
\nabla u_\ep^+ = {U_0}' \displaystyle{\frac{\nabla d}{\ep}},\vsp \\
\Delta u_\ep^+= {U_0}''\displaystyle{\frac{|\nabla d|^2}{\ep ^2}}
+ {U_0}'\displaystyle{\frac{\Delta d}{\ep}},
\end{array}
$$
where the function $U_0$, as well as its derivatives, is taken at
the point $\big(d (x,t)-\ep p(t)\big)/ \ep $. We also use
expansions of the reaction terms:
$$
\begin{array}{ll}
f(u_\ep^+)=f(U_0)+qf '(U_0) + \di{\frac 12}q^2f ''(\theta), \vsp \\
g(u_\ep^+)=g(U_0)+q g'(\omega),
\end{array}
$$
where $\theta(x,t)$ and $\omega(x,t)$ are some functions
satisfying $U_0<\theta<u_\ep^+,\;U_0<\omega<u_\ep^+$. Combining
the above expressions with equation \eqref{eq-phi-chemo} and the
fact that $\sqrt 2 \alpha {U_0}'+g(U_0)\equiv 0$, we obtain
$$
L_{v^\ep}[u_\ep^+]=E_1+\cdots+E_5,
$$
where:\vsp\\
$E_1=- \edeux q[f '(U_0)
+\frac 12 q f ''(\theta)]-{U_0}'p_t+q_t$,\vsp\\
$E_2=\di{\frac{{U_0}''}{\ep^2}}(1-|\n d|^2)$,\vsp\\
$E_3=\displaystyle{\frac {{U_0}'}{\ep}}(d_t-\Delta d +
\n d\cdot \n \chi (v^0)+\sqrt 2 \alpha)$,\vsp \\
$E_4=-\displaystyle{\frac{1}{\ep}}
qg'(\omega)$,\vsp \\
$E_5=\di{\frac{{U_0}'}{\ep}}\n d\cdot \n (\chi(v^\ep)-\chi(v^0))+
u_ \ep ^+\Delta \chi(v^\ep).$ \\

In order to estimate the terms above, we first present some useful
inequalities. As $f'(0)=f'(1)=-1/2$, we can find strictly positive
constants $b$ and $m$ such that
\begin{equation}\label{bords-chemo}
\textrm{ if }\quad U_0(z) \in [0,b]\cup [1-b,1] \quad \quad
\textrm{ then }\quad f'(U_0(z))\leq -m.
\end{equation}
On the other hand, since the region $\{z\in\R\,|\,U_0(z)\in
[b,1-b] \,\}$ is compact and since ${U_0}'<0$ on $\R$, there
exists a constant $a_1>0$ such that
\begin{equation}\label{milieu-chemo}
\textrm{ if }\quad U_0(z)\in [b,1-b] \quad \textrm{ then } \quad
{U_0}'(z) \leq - a_1.
\end{equation}
We then define
\begin{equation}\label{F-chemo}
F=\sup_{-1\leq z\leq 2} |f(z)|+|f'(z)|+|f ''(z)|,
\end{equation}
\begin{equation}\label{beta-chemo}
\beta = \frac{m}{4},
\end{equation}
and choose $\sigma$ that satisfies
\begin{equation}\label{sigma-chemo}
0< \sigma \leq \min\,(\sigma_0,\sigma_1,\sigma_2),
\end{equation}
where
\[
\sigma _0:=\frac{a_1}{m+F},\quad \sigma _1:=\frac{1}{\beta
+1},\quad \sigma _2:=\frac{4\beta }{F(\beta+1)}.
\]
Hence, combining \eqref{bords-chemo} and \eqref{milieu-chemo}, we
obtain, using that $\sigma \leq \sigma _0$,
\begin{equation}\label{U0-f-chemo}
-{U_0}'(z)-\sigma f'(U_0(z))\geq 4 \sigma \beta \qquad \hbox{for}
\ \ -\infty<z<\infty.
\end{equation}

Now let $K>1$ be arbitrary. In what follows we will show that
$L_{v^\ep}[u _\ep ^+] \geq 0$ provided that the constants $\ep_0$
and $L$ are appropriately chosen. From now on, we suppose that the
following inequality is satisfied:
\begin{equation}\label{ep0M-chemo}
\ep _0^2 Le^{LT} \leq 1\, .
\end{equation}
Then, given any $\ep\in(0,\ep_0)$, since $\sigma \leq \sigma _1$,
we have $0\leq q(t)\leq 1$, hence, recalling that $0<U_0<1$,
\begin{equation}\label{uep-pm-chemo}
-1\leq u_\ep^\pm(x,t) \leq 2\, .
\end{equation}

\vskip 8 pt \noindent{ \bf We first estimate the term $E_1$.} A
direct computation gives
$$
E_1=\frac{\beta}{\ep^2}\,\EB(I-\sigma\beta)+Le^{Lt}(I+\ep^2\sigma
L),
$$
where
$$
I=-{U_0}'-\sigma f '(U_0)-\frac {\sigma^2}2 f
''(\theta)(\beta\EB+\ep^2 Le^{Lt}).
$$
In virtue of \eqref{U0-f-chemo} and \eqref{uep-pm-chemo}, we
obtain
\[
I\geq 4 \sigma \beta-\frac {\sigma^2}{2} F(\beta+\ep^2 Le^{LT}).
\]
Then, in view of \eqref{ep0M-chemo}, using that $\sigma \leq
\sigma _2$, we have
\[
I\geq 2\sigma\beta.
\]
Consequently, the following inequality holds.
$$
E_1\geq \frac{\sigma\beta^2}{\ep^2}\EB + 2\sigma\beta L
e^{Lt}=:\frac {C_1}{\ep^2}\EB + {C_1}'L e^{Lt}.
$$

\vskip 8 pt \noindent{\bf The term $E_2$.} First, in the points
where where $|d(x,t)|<d_0$, we have that $|\n d|=1$ so that
$E_2=0$. Next we consider the points where $|d(x,t)|\geq d_0.$ We
deduce from Lemma \ref{est-phi-chemo} that:
$$
\begin{array}{ll}
|E_2|&\leq \di{\frac{C}{\ep^2}}(1+\Vert \n d\Vert _\infty
^2)e^{-\lambda|d+\ep p|/ \ep}\vsp\\
&\leq \di{\frac{C}{\ep^2}}(1+\Vert \n d\Vert _\infty
^2)e^{-\lambda(d_0 / \ep-|p|)}.
\end{array}
$$
In view of the definition of $p$ in \eqref{defpetq-chemo}, we have
that $0<K-1 \leq p \leq e^{LT} +K$, and suppose from now that the
following assumption holds:
\begin{equation}\label{cond3-chemo}
e^{LT}+K \leq \frac{d_0}{2\ep_0}.
\end{equation}
Then $\displaystyle{\frac{d_0}{\ep}}-|p|\geq
\displaystyle{\frac{d_0}{2\ep}}$, so that
$$
\begin{array}{ll}
|E_2|&\leq \di{\frac{C}{\ep^2}}(1+\Vert \n d\Vert _\infty
^2)e^{-\lambda d_0 / (2\ep)}\vsp\\
&\leq C_2:=\di{\frac{16C}{(e\lambda d_0)^2}}(1+\Vert \n d\Vert
_\infty ^2).
\end{array}
$$

\vskip 8 pt \noindent{\bf Next we consider the term $E_3$.} We
recall that
$$
d_t-\Delta d +\n d \cdot \n\chi(v^0)+\sqrt 2 \alpha=0 \quad\
\textrm{on} \ \; \Gamma _t=\{x \in \om, d(x,t)=0\}.
$$
Since $v^0$ is of class $C^{1+\vartheta ',\frac{1+\vartheta '}2}$,
for any $\vartheta ' \in (0,1)$, and since the interface $\gam _t$
is of class $C^{2+\vartheta,\frac{2+\vartheta}2}$, the functions
$\n d$, $\Delta d$, $d_t$ and $\n \chi(v^0)$ are Lipschitz
continuous near $\gam _t$. It then follows, from the mean value
theorem applied separately on both sides of $\Gamma_t$, that there
exists $N_0>0$ such that:
$$
|(d_t-\Delta d +\n d \cdot \n\chi(v^0)+\sqrt 2 \alpha)(x,t)|\leq
N_0|d(x,t)| \quad\ \textrm{for all} \ \; (x,t) \in Q_T.
$$
Applying  Lemma \ref{est-phi-chemo} we deduce that
$$
\begin{array}{lll}
|E_3|&\leq N_0C\displaystyle{\frac
{|d(x,t)|}{\ep}}e^{-\lambda| d(x,t)/ \ep +p(t)|}\vsp \\
&\leq N_0C \max_{y \in \mathbb{R}}|y|e^{-
\lambda|y +p(t)|}\vsp \\
&\leq N_0C\max \big(|p(t)|,\di{\frac 1 \lambda}\big)\vsp\\
&\leq N_0C\big(|p(t)|+\di{\frac 1 \lambda}\big).
\end{array}
$$
Taking the expression of $p$ into account, we see that $|p(t)|\leq
e^{Lt}+K$, which implies
$$
|E_3|\leq C_3(e^{Lt}+K)+{C_3}',
$$
where $C_3:=N_0C$ and ${C_3}':=N_0C/\lambda$.

\vskip 8 pt \noindent{\bf The term $E_4$.} We set $G_1:=\Vert g'
\Vert _{L^\infty (-1,2)}$ and, substituting the expression for
$q$, obtain that
$$
\begin{array}{ll}
|E_4| &\leq \sigma G_1\Big(\di{\frac \beta \ep} e^{-\beta
t/\ep^2}+\ep Le^{Lt}\Big)\vsp \\
&\leq\di{ \frac {C_4}{\ep}} e^{-\beta t/\ep^2}+{C_4}'\ep Le^{Lt}.
\end{array}
$$

\vskip 8 pt \noindent{\bf We continue with the term $E_5$.} This
term requires a more delicate analysis.  We need a precise
estimate of $v^\ep-v^0$. We recall that $v^0$ satisfies $-\Delta
v^0 +\gamma v^0=u^0$, with $u^0$ a step function discontinuous
when crossing the interface.

\begin{lem}\label{green-chemo}There exists a positive constant $C_G$ such that, for all $(x,t)\in
\Q$,
\begin{equation}\label{elliptic-estimate-chemo}
\Big(|v^\ep|+|\n v^\ep|+|\Delta v^\ep|\Big) (x,t) \leq C_G ,
\end{equation}
\begin{equation}\label{est-green-chemo}
\Big(|v^\ep-v^0|+| \n d \cdot \n (v^\ep-v^0)|\Big)(x,t) \leq
C_G(\ep p(t) +q(t)).
\end{equation}
\end{lem}

\vskip 8 pt \noindent We postpone the proof of this lemma and
pursue the proof of Lemma \ref{fix-chemo}. Using the smoothness of
$\chi$ and \eqref{elliptic-estimate-chemo}, we obtain a uniform
bound ${C_G}'$ for $\Delta \chi(v^\ep)$. Moreover, we write
\begin{equation}\label{chilogak-chemo}
\n d\cdot \n \big(\chi (v^\ep)-\chi (v^0)\big)={\chi} '(v^\ep)\n
d\cdot \n (v^\ep -v^0)+\big(\chi '(v^\ep)-\chi '(v^0)\big) \n
d\cdot \n v^0.
\end{equation}
Since $v^0$ is of class $C^{1+\vartheta ',\frac{1+\vartheta '}2}$,
for any $\vartheta ' \in (0,1)$, there exists a constant, which we
denote again by $C_G$, such that
$$
\Vert v^0\Vert _{L^\infty(Q_T)}+ \Vert \n v ^0 \Vert
_{L^\infty(Q_T)} \leq C_G,
$$
which, combined with \eqref{chilogak-chemo}, yields
$$
|\n d\cdot \n\big (\chi (v^\ep)-\chi (v^0)\big)|\leq \Vert \chi '
\Vert _\infty |\n d \cdot \n (v^\ep -v^0)|\\
+|v^\ep -v^0|\ \Vert \chi '' \Vert _\infty \Vert \n d\Vert _\infty
C_G,
$$
where the $L^\infty$-norms of $\chi '$ and $\chi ''$ are
considered on the interval $(-C_G,C_G)$. It follows from the above
inequality and \eqref{est-green-chemo} that there exists a
constant ${C_G}''$ such that, for all $(x,t) \in Q_T$,
$$
|\n d\cdot \n\big (\chi (v^\ep)-\chi (v^0)\big)|(x,t) \leq
{C_G}''(\ep p(t) +q(t)).
$$
Hence, using \eqref{uep-pm-chemo} and the above estimates, we
obtain,
$$
|E_5|\leq \frac{C}{\ep}{C_G}''(\ep p(t)+q(t))+2{C_G}'.
$$
Then, substituting the expressions for $p$ and $q$, we easily
obtain positive constants $C_5$, ${C_5}'$ and ${C_5}''$ such that
$$
|E_5|\leq C_5+ \frac {{C_5}'}{\ep}e^{-\beta t/
\ep^2}+{C_5}''(1+\ep  L)e^{Lt}.
$$

\vskip 8 pt \noindent {\bf Completion of the proof.} Collecting
the above estimates of $E_1$---$E_5$ yields
$$
L_{v^\ep}[u_\ep^+]\geq \frac{C_1-\ep
C_4-\ep{C_5}'}{\ep^2}e^{-\beta t/\ep ^2}+
\big(L({C_1}'-\ep{C_4}'-\ep {C_5}'')-C_3-{C_5}''\big)e^{Lt}-C_7,
$$
where $ C_7:=C_2+KC_3+{C_3}'+C_5$. Now, we set
$$
L:=\frac 1 T\ln \frac {d_0}{4\ep _0},
$$
which, for $\ep_ 0$ small enough, validates assumptions
\eqref{ep0M-chemo} and \eqref{cond3-chemo}. If $\ep_0$ is chosen
sufficiently small (i.e. $L$ large enough), $C_1/\ep
^2-(C_4+{C_5}')/ \ep$ is positive, ${C_1}'-\ep{C_4}'-\ep
{C_5}''\geq \frac 1 2 {C_1}'$, and
$$
\begin{array}{lll}
L_{v^\ep}[u_\ep ^+]&\geq \big[\frac 1 2 L{C_1}'-C_3-{C_5}''\big]e^{Lt}-C_7\vsp \\
&\geq  \frac 1 4 L{C_1}' -C_7 \vsp \\
&\geq 0.
\end{array}
$$
The proof of Lemma \ref{fix-chemo} is now completed, with the
choice of the constants $\beta, \sigma$ as in \eqref{beta-chemo},
\eqref{sigma-chemo}.\qed

\subsection{Proof of Lemma \ref{green-chemo}}

Lemma \ref{green-chemo} is inspired by Lemma 4.9 in \cite{BHLM3}.
Since our pair of sub- and super-solutions is different from the
one in \cite{BHLM3}, we need to perform some minor changes. First
we give a useful estimate on ``shifted $U_0$".

\begin{lem}\label{est-phi2-chemo}
For all $a\in\mathbb{R}$, all $z\in\mathbb{R}$, we have
$$|U_0(z+a)-\chi_{]-\infty,0]}(z)|\leq
Ce^{-\lambda|z+a|}+\chi_{[-a,a]}(z)
$$
\end{lem}
{\noindent \bf Proof.} Let us give the proof for $a>0$. We
distinguish three cases and use the estimates of Lemma
\ref{est-phi-chemo}. For $z\leq -a$, we have $|U_0(z+a)-1|\leq
Ce^{-\lambda|z+a|}$. For $-a < z \leq 0$, we have
$|U_0(z+a)-1|\leq |U_0(z+a)|+1 \leq Ce^{-\lambda|z+a|}+1$. For
$z>0$, we have $|U_0(z+a)|\leq Ce^{-\lambda|z+a|}$. We proceed in
the same way for $a<0$. \qed

\vskip 8pt We turn to the proof of Lemma \ref{green-chemo}. First,
we recall that $v^\ep$ is such that \eqref{v-coincee-def-chemo}
holds; hence, in view of \eqref{uep-pm-chemo}, the estimate
\eqref{elliptic-estimate-chemo} is a direct consequence of the
standard theory of elliptic equations. Next we prove
\eqref{est-green-chemo}. The function $w=w^\ep:=v^\ep-v^0$ is
solution of
\begin{equation}\label{eq-w-chemo}
\begin{cases}
-\Delta w + \gamma w =h &\textrm{ on } \Q, \vsp\\
\di{\frac{\partial w}{\partial \nu}}=0 &\textrm{ on } \partial
\Omega \times (0,T),
\end{cases}
\end{equation}
with $u_\ep ^- -u^0 \leq h=h^\ep \leq u_\ep ^+ -u^0$, where $u^0$
is the step function defined by $u^0(x,t)=\chi_{\{d(x,t)\leq
0\}}$. The key idea of the proof is the fact that $h$ is
exponentially small with respect to $\ep$, except possibly in a
thin neighborhood of $\Gamma _t$ of width of order $\ep p(t)$.
More precisely, from the definitions of $u_ \ep^{\pm}$ in
\eqref{sub-chemo} and from the above lemma for $z=d(x,t)/\ep$ and
$a=\pm p(t)$, we deduce that
\begin{equation}\label{est-h-chemo}
|h(x,t)| \leq C (e^{-\lambda | d(x,t)/ \ep +p(t)|}+e^{-\lambda |
d(x,t)/ \ep-p(t)|})+ \chi _{\{|d(x,t)| \leq \ep p(t)\}}+q(t).
\end{equation}
By linearity, we successively consider equation \eqref{eq-w-chemo}
with the various terms appearing in the right-hand side of
\eqref{est-h-chemo}. By the standard elliptic estimates, the
solution $w$ of \eqref{eq-w-chemo} satisfies
\begin{equation}\label{ineq-w-chemo}
|w(x,t)|+|\n w(x,t)| \leq C' \sup _{y \in \Omega} |h(y,t)|,
\end{equation}
which gives the term $C_G q(t)$ that appears in the right-hand
side of inequality \eqref{est-green-chemo} for $h(y,t)=q(t)$. We
now suppose that the function $h$ satisfies one of the three
following assumptions:
$$
\begin{array}{ll}
(H_1) \quad \quad &|h(y,t)| \leq \di{\chi} _{\di{\{}|d(y,t)| \leq  \ep p(t)\di{\}}} \vsp\\
(H_2 ^\pm) \quad \quad &|h(y,t)| \leq  \exp
\Big(-\lambda|\di{\frac{d(y,t)}\ep }\pm  p(t)|\Big),
\end{array}
$$
and write
$$
h(y,t)=h(y,t)\chi_{\{|d(y,t)|\leq d_0\}}+h(y,t)\chi_{\{|d(y,t)|>
d_0\}}.
$$

We first consider the term $h(y,t)\chi_{\{|d(y,t)|> d_0\}}$. In
virtue of \eqref{cond3-chemo}, we have
\begin{equation}\label{p-coincee-chemo}
0<K-1\leq p(t)\leq d_0 /2\ep_0.
\end{equation}
Under assumption $(H_1)$, it follows that $h$ is supported in
$\{|d(y,t)|\leq d_0/2\}$, which implies $h(y,t)\chi_{\{|d(y,t)|>
d_0\}}=0$. Moreover, under assumption $(H_2 ^\pm)$, using again
\eqref{p-coincee-chemo},
$$
\begin{array}{llll}
|h(y,t)|\chi_{\{|d(y,t)|> d_0\}} &\leq \exp \big [-\lambda(d_0
/\ep-p(t))\big]\vsp\\
& \leq \exp (-\lambda d_0 / 2\ep)  \vsp\\
& \leq \di{\frac {2}{\lambda d_0 e}}\; \ep\vsp  \\&\leq \di{\frac
{2}{\lambda d_0 e}}\frac {1}{K-1}\ep p(t).
\end{array}
$$
Thus, under either of the assumptions $(H_1)$ or $(H_2^\pm)$, the
estimate \eqref{est-green-chemo} --- for the term
$h(y,t)\chi_{\{|d(y,t)|>d_0\}}$ --- directly follows from
inequality \eqref{ineq-w-chemo}.

From now on, we assume that $h$ is supported in $\{|d(y,t)| \leq
d_0\}$. We have that
$$
w(x,t)=\int _{|d(y,t)|\leq d_0} G(x,y)h(y,t)dy,
$$
and
$$\n d(x,t) \cdot \n w(x,t)= \int _{|d(y,t)|\leq d_0}(\n _x
G(x,y)\cdot \n d (x,t))h(y,t)dy,
$$
where $G$ is the Green's function associated to the homogeneous
Neumann boundary value problem on $\om$ for the operator $-\Delta+
\ga$. More precisely, $G(x,y)=g_\gamma(|x-y|)+H_\gamma(x,y)$,
where $g_\gamma(|x-y|)$ is the Green's function associated to the
operator $-\Delta +\gamma$ on $\R ^N$ and where $H_\gamma (x,y)$
is smooth for $x$ and $y$ far away from $\partial \om$. It is
known that $g_\gamma$ is the Bessel function defined by
$$
g_\gamma(r)=c_N\int _0 ^{+\infty} e^{-\frac{r^2}{2s}}e^{-\gamma
\frac s 2}s^{\frac{-N+2}{2}}\frac{ds}{s},
$$
with $c_N>0$ a normalization constant. We use the following
estimates (see \cite{BHLM3}):
\begin{equation}\label{est-G-chemo}
|G(x,y)|\leq \begin{cases}\di{\frac{C}{|y-x|^{N-2}}} \quad
&\text{ for } N\geq 3 \vsp \\
C | \ln |y-x|| \quad &\text{ for } N=2,\end{cases}
\end{equation}
\begin{equation}\label{est-nabla-G-chemo}
|\n _xG(x,y)\cdot \n d(x,t)|\leq
\di{\frac{C|d(y,t)-d(x,t)|}{|y-x|^N}}+\di{\frac{C}{|y-x|^{N-2}}}\quad
\text{ for } N\geq 2.
\end{equation}
This last inequality follows from
$$
|\n _xG(x,y)\cdot \n d(x,t)|\leq \di{\frac{C|\n d(x,t)\cdot
(y-x)|}{|y-x|^N}},
$$
and from $d(y,t)-d(x,t)=\n d(x,t) \cdot (y-x)+O(|y-x|^2)$. Now,
under respectively assumptions $(H_1)$, $(H_2 ^\pm)$, we define a
function $\tilde h=\tilde h^\ep$ on $\R \times [0,T]$,
respectively by
\begin{equation}\label{tilde-h-chemo}
\tilde h (r,t):= \begin{cases}\chi_{\di{\{}|r|\leq \ep p(t)\di{\}}}  \vsp \\
\exp \Big(-\lambda|\di{\frac r \ep} \pm p(t)|\Big).\end{cases}
\end{equation}
Note that $|h(y,t)| \leq \tilde h\big(d(y,t),t\big)$. Moreover,
using \eqref{p-coincee-chemo}, straightforward computations show
that, under either of the assumptions $(H_1)$ or $(H_2 ^\pm)$,
there exists $\tilde C>0$ such that
\begin{equation}\label{int-h-chemo}
0\leq \int _{-d_0}^{d_0} \tilde h (r,t)dr \leq \tilde C \ep p(t),
\end{equation}
which is an analogue of $(4.20)$ in \cite{BHLM3}. The end of the
proof is now identical to that of Lemma 4.10 in \cite{BHLM3}. We
omit the details and refer to this article. \qed

\section{Proof of the main results }\label{s:proof-chemo}

In this section, we prove our main results by fitting the two
pairs of sub- and super-solutions, constructed for the study of
the generation and the motion of interface, into each other.

\subsection{Proof of Theorem \ref{width-chemo}}\label{ss:proof-w-chemo}

Let $\eta \in (0,1/4)$ be arbitrary. Choose $\beta$ and $\sigma$
that satisfy \eqref{beta-chemo}, \eqref{sigma-chemo} and
\begin{equation}\label{eta-chemo}
\sigma \beta \leq \frac \eta 3.
\end{equation}
By the generation of interface Theorem \ref{th-gen-chemo}, there
exist positive constants $\ep_0$ and $M_0$ such that
\eqref{part1-chemo}, \eqref{part2-chemo} and \eqref{part3-chemo}
hold with the constant $\eta$ replaced by $ \sigma \beta /2$.
Since $\n u_0 \cdot n \neq 0$ everywhere on the initial interface
$\Gamma _0=\{x\in\om, \; u_0(x)=1/2 \}$ and since $\Gamma _0$ is a
compact hypersurface, we can find a positive constant $M_1$ such
that
\begin{equation}\label{corres-chemo}
\begin{array}{ll}\text { if } \quad d_0 (x) \geq \quad M_1 \ep &
\text { then } \quad
u_0(x) \leq 1/2 -M_0 \ep,\vspace{3pt}\\
\text { if } \quad d_0 (x) \leq -M_1 \ep & \text { then } \quad
u_0(x) \geq 1/2+M_0 \ep.
\end{array}
\end{equation}
Here $d_0(x):=d(x,0)$ denotes the cut-off signed distance function
associated with the hypersurface $\Gamma_0$. Now we define
functions $H^+(x), H^-(x)$ by
\[
\begin{array}{l}
H^+(x)=\left\{
\begin{array}{ll}
1+\frac 12\sigma\beta\quad\ &\hbox{if}\ \ d_0(x)< M_1\ep\vsp \\
\;\;\;\;\;\;\frac 12 \sigma\beta\quad\ &\hbox{if}\ \ d_0(x)\geq
M_1\ep,
\end{array}\right.
\vsp\\
H^-(x)=\left\{
\begin{array}{ll}
1-\frac 12\sigma\beta\quad\ &\hbox{if}\ \ d_0(x)\leq \;-M_1\ep\vsp \\
\;\;\;-\frac 12\sigma\beta\quad\ &\hbox{if}\ \ d_0(x) > \;-M_1\ep.
\end{array}\right.
\end{array}
\]
Then from the above observation we see that
\begin{equation}\label{H-u-chemo}
H^-(x) \,\leq\, u^\ep(x,\mu^{-1} \ep^2|\ln \ep|) \,\leq\,
H^+(x)\qquad \hbox{for}\ \ x\in\Omega.
\end{equation}

Next we fix a sufficiently large constant $K>1$ such that
\begin{equation}\label{K-chemo}
U_0(M_1-K) \geq 1-\frac {\sigma \beta}{3} \quad \text { and }
\quad U_0(-M_1+K) \leq \frac {\sigma \beta}{3}.
\end{equation}
For this $K$, we choose $\ep _0$ and $L$  as in Lemma
\ref{fix-chemo}. We claim that
\begin{equation}\label{uep-H-chemo}
u_\ep^-(x,0)\leq H^-(x),\quad\ H^+(x)\leq u_\ep^+(x,0) \qquad
\hbox{for} \ \ x\in\Omega.
\end{equation}
We only prove the former inequality, as the proof of the latter is
virtually the same. Then it amounts to showing that
\begin{equation}\label{c3-chemo}
u_\ep^-(x,0)=U_0\big(\frac {d_0(x)}{\ep}+K\big)-\sigma (\beta+\ep
^2 L) \;\leq\; H^-(x).
\end{equation}
In the range where $d_0(x) >- M_1 \ep$, the second inequality in
\eqref{K-chemo} and the fact that $U_0$ is a decreasing function
imply
\[
U_0\big(\frac {d_0(x)}{\ep}+K\big)-\sigma (\beta+\ep ^2 L)\leq
\frac{1}{3}\sigma\beta -\sigma\beta\;\leq\;H^-(x).
\]
On the other hand, in the range where $d_0(x) \leq- M_1 \ep$, we
have
\[
U_0\big(\frac {d_0(x)}{\ep}+K\big)-\sigma (\beta+\ep ^2 L)\leq
\;1-\sigma\beta\;\leq\;H^-(x).
\]
This proves \eqref{c3-chemo}, hence \eqref{uep-H-chemo} is
established.

Combining \eqref{H-u-chemo} and \eqref{uep-H-chemo}, we obtain
$$
u_\ep^-(x,0)\leq u^\ep(x,\mu ^{-1} \ep ^2|\ln \ep|) \leq
u_\ep^+(x,0).
$$
Since $u_\ep^-$ and $u_\ep^+$ are sub- and super-solutions for
Problem $\Pe$ thanks to Lemma \ref{fix-chemo}, the comparison
principle yields
\begin{equation}\label{ok-chemo}
u_\ep^-(x,t) \leq u^\ep (x,t+t^\ep) \leq u_\ep^+(x,t) \quad \text
{ for } \quad 0 \leq t \leq T-t^\ep,
\end{equation}
where $t^\ep=\mu ^{-1} \ep ^2|\ln \ep|$. Note that, in view of
\eqref{sub2-chemo}, this is enough to prove Corollary
\ref{total-chemo}. Now let $C$ be a positive constant such that
\begin{equation}\label{C-chemo}
U_0(-C+e^{LT}+K) \geq 1-\frac \eta 2 \quad \text { and } \quad
U_0(C-e^{LT}-K) \leq \frac \eta 2.
\end{equation}
One then easily checks, using successively \eqref{ok-chemo},
\eqref{sub-chemo}, \eqref{C-chemo} and \eqref{eta-chemo}, that,
for $\ep _0$ small enough, for $0 \leq t \leq T-t^\ep$, we have
\begin{equation}\label{correspon-chemo}
\begin{array}{ll}\text { if } \quad d(x,t) \geq \quad C \ep &
\text { then } \quad
u^\ep (x,t+t^\ep) \leq \eta,\vspace{3pt}\\
\text { if } \quad d(x,t) \leq -C \ep & \text { then } \quad u^\ep
(x,t+t^\ep) \geq 1-\eta,
\end{array}
\end{equation}
and
$$
u^\ep (x,t+t^\ep) \in [-\eta,1+\eta],
$$
which completes the proof of Theorem \ref{width-chemo}.\qed

\subsection{Proof of Theorem \ref{error-chemo}}

In the case where $ \mu ^{-1} \ep ^2 |\ln \ep|\leq t \leq T$, the
assertion of the theorem is a direct consequence of Theorem
\ref{width-chemo}.  All we have to consider is the case where
$0\leq t \leq \mu ^{-1} \ep ^2 |\ln \ep|$. We shall use the sub-
and super-solutions constructed for the study of the generation of
interface in Section \ref{generation-chemo}. To that purpose, we
first prove the following lemma concerning $Y(\tau,\xi;\de)$, the
solution of the ordinary differential equation \eqref{ode-chemo},
in the initial time interval.

\begin{lem}\label{est-Y-debut-chemo}
There exist constants $C_8 >0$ and $\ep _0>0$ such that, for all
$\ep \in(0,\ep _0)$,
\begin{equation}\label{debut-chemo}
\begin{array}{ll}\text { if } \quad \xi \geq 1/2+C_8 \ep &
\text{then} \quad Y(\tau,\xi;\pm \ep  G)> 1/2\quad \hbox{for}
\quad 0\leq \tau \leq \mu ^{-1}|\ln \ep|,
\vspace{3pt}\\
\text { if } \quad \xi \leq 1/2-C_8 \ep & \text {then} \quad
Y(\tau,\xi;\pm \ep G)<1/2 \quad \hbox{for} \quad 0\leq \tau \leq
\mu ^{-1}|\ln \ep|.
\end{array}
\end{equation}
\end{lem}

{\noindent \bf Proof.} We only prove the first inequality. Assume
$\xi \geq 1/2+C_8 \ep$. By \eqref{h-chemo}, for $C_8 \geq C G$, we
have that $\xi\geq 1/2+C_8 \ep \geq a(\pm\ep G)$. It then follows
from \eqref{est-Y-1-chemo} that
$$
\begin{array}{llll}Y(\tau,\xi;\pm \ep G)&\geq a(\pm \ep G)+C_1 e^{\mu(\pm \ep
G)\tau}(1/2+C_8 \ep-a(\pm \ep G)) \vsp \\ &\geq 1/2-CG \ep+C_1(-CG
\ep+C_8 \ep)\vsp \\ &\geq 1/2+\ep(C_1C_8-CG(C_1 +1))\vsp
\\ &>1/2,
\end{array}
$$
provided that $C_8$ is sufficiently large. \qed

\vskip 8pt Now we turn to the proof of Theorem \ref{error-chemo}.
We first claim that there exists a positive constant $M_2$ such
that for all $t \in [0,\mu ^{-1} \ep ^2 |\ln \ep|]$,
\begin{equation}\label{claim1-chemo}
\Gamma _t ^\ep \subset \mathcal N _{M_2\ep} (\Gamma _0).
\end{equation}
To see this, we choose ${M_0}'$ large enough, so that ${M_0}' \geq
C_8+2 C_6$, where $C_6$ is as in Lemma \ref{w-chemo}. As is done
for \eqref{corres-chemo}, there is a positive constant $M_2$ such
that
\begin{equation}\label{corres2-chemo}
\begin{array}{ll}\text { if } \quad d_0 (x) \geq \quad M_2 \ep &
\text { then } \quad
u_0(x) \leq 1/2 -{M_0}' \ep,\vspace{3pt}\\
\text { if } \quad d_0 (x) \leq -M_2 \ep & \text { then } \quad
u_0(x) \geq 1/2 +{M_0}' \ep.
\end{array}
\end{equation}
In view of this last condition, we see that, if $\ep _0$ is small
enough, if $d_0(x)\geq M_2 \ep$, then for all $0 \leq t \leq
\mu^{-1} \ep ^2 | \ln \ep|$,
$$
\begin{array}{llll}u_0(x)+\ep ^2 r(\ep G,\di{\frac
{t}{\ep ^2}})&\leq 1/2-{M_0}' \ep+\ep^2 C_6\big[e^{\mu(\ep
G)|\ln \ep|/ \mu}-1\big]\vsp \\
&\leq 1/2+\ep\big[-{M_0}'+C_6 \ep ^{(\mu-\mu(\pm\ep G))/\mu}-\ep
C_6\big]\vsp \\
&\leq 1/2+\ep(-{M_0}' +2 C_6) \quad\quad  \leftarrow \text {thanks
to \eqref{point-chemo}}\vsp
\\ &\leq 1/2-C_8 \ep.
\end{array}
$$
This inequality and Lemma \ref{est-Y-debut-chemo} imply
$w_\ep^+(x,t)<1/2$, where $w_\ep^+$ is the sub-solution defined in
\eqref{w+--general-chemo}.  Consequently, by
\eqref{coincee-chemo},
\[
u^\ep(x,t)<1/2\qquad \hbox{if} \ \ d_0(x)\geq M_2 \ep.
\]
In the case where $d_0(x)\leq -M_2\ep$, similar arguments lead to
$u^\ep (x,t)>1/2$. This completes the proof of
\eqref{claim1-chemo}. Note that we have proved that, for all $0
\leq t \leq \mu ^{-1} \ep ^2 |\ln \ep|$,
\begin{equation}\label{rost1-chemo}
\begin{array}{ll}
u^\ep(x,t) >1/2 \quad \text {if} \quad x \in \Omega _0 ^{(1)}
\setminus
\mathcal N _{M_2\ep}(\Gamma _0),\vspace{3pt}\\
u^\ep(x,t) <1/2 \quad \text {if} \quad x \in \Omega _0 ^{(0)}
\setminus \mathcal N _{M_2\ep}(\Gamma _0).
\end{array}
\end{equation}

Next, since $\Gamma _t$ depends on $t$ smoothly, there is a
constant $\tilde C >0$ such that, for all $t \in [0,\mu ^{-1} \ep
^2 |\ln \ep|]$,
\begin{equation}\label{claim2-chemo}
\Gamma _0 \subset \mathcal N _{\tilde C\ep ^2|\ln \ep|} (\Gamma
_t),
\end{equation}
and
\begin{equation}\label{rost2-chemo}
\begin{array}{ll}
\Omega _t ^{(1)} \setminus \mathcal N _{\tilde C \ep}(\Gamma _t)
\subset
\Omega _0 ^{(1)} \setminus \mathcal N _{M_2\ep}(\Gamma _0),\vspace{3pt}\\
\Omega _t ^{(0)} \setminus \mathcal N _{\tilde C \ep}(\Gamma _t)
\subset \Omega _0 ^{(0)} \setminus \mathcal N _{M_2\ep}(\Gamma
_0).
\end{array}
\end{equation}

As a consequence of \eqref{claim1-chemo} and \eqref{claim2-chemo}
we get
$$
\Gamma _t ^\ep \subset \mathcal N _{M_2 \ep+\tilde C \ep ^2|\ln
\ep|} (\Gamma _t) \subset \mathcal N _{C \ep} (\Gamma _t),
$$
which completes the proof of Theorem \ref{error-chemo}.\qed

\vskip 8 pt {\bf \noindent Proof of Corollary
\ref{total-2-chemo}.} In view of Theorem \ref{error-chemo} and the
definition of the Hausdorff distance, to prove this corollary we
only need to show the reverse inclusion, that is
\begin{equation}\label{obj-chemo}
\Gamma _t \subset \mathcal N _{C' \ep} (\Gamma _t ^\ep)\quad \text
{ for } \quad 0 \leq t \leq T,
\end{equation}
for some constant $C'>0$. To that purpose let $C'$ be a constant
satisfying $C'>\max (\tilde C,C)$, where $C$ is as in Theorem
\ref{width-chemo} and $\tilde C$ as in \eqref{rost2-chemo}. Choose
$t \in [0,T]$, $x_0 \in \Gamma _t$ arbitrarily and, $n$ being the
Euclidian normal vector exterior to $\Gamma _t$ at point $x_0$,
define a pair of points:
$$
x^{(0)}:=x_0+C'\ep n \quad \text{and}\quad  x^{(1)}:=x_0-C' \ep n.
$$
Since $C'>C$ and since the curvature of $\Gamma _t$ is uniformly
bounded as $t$ varies over $[0,T]$, we see that, if $\ep _0$ is
sufficiently small,
$$
x^{(0)} \in \Omega ^{(0)} _t \setminus \mathcal N _{C\ep}(\Gamma
_t) \quad \text {and} \quad  x^{(1)} \in \Omega ^{(1)} _t
\setminus \mathcal N _{C\ep}(\Gamma _t).
$$
Therefore, if $t\in [\mu ^{-1} \ep ^2|\ln \ep|,T]$, then, by
Theorem \ref{width-chemo}, we have
\begin{equation}\label{rostand-chemo}
u^\ep(x^{(0)},t)<1/2<u^\ep(x^{(1)},t).
\end{equation}
On the other hand, if $t\in[0,\mu ^{-1} \ep ^2|\ln \ep|]$, then
from \eqref{rost1-chemo}, \eqref{rost2-chemo} and the fact that
$C'>\tilde C$, we again obtain \eqref{rostand-chemo}. Thus
\eqref{rostand-chemo} holds for all $t\in[0,T]$. Now, by the mean
value theorem, we see that, for each $t\in[0,T]$, there exists a
point $\tilde x$ such that
$$
\tilde x \in [x^{(0)},x^{(1)}] \quad \text{and} \quad u^\ep(\tilde
x,t)=1/2.
$$
This implies $\tilde x \in \Gamma _t ^\ep$. Furthermore we have
$|x_0-\tilde x|\leq C' \ep$, since $\tilde x$ lies on the line
segment $[x^{(0)},x^{(1)}]$. This proves \eqref{obj-chemo}. \qed


\begin{thebibliography}{ABCD}

\bibitem{AHM} M.~Alfaro, D.~Hilhorst and H.~Matano,
{\it The singular limit of the Allen-Cahn equation and the
FitzHugh-Nagumo system}, in preparation.

\bibitem{AC} S.~Allen and J.~Cahn,
{\it A microscopic theory for antiphase boundary motion and its
application to antiphase domain coarsening}, Acta Metall. {\bf 27}
(1979), 1084--1095.

\bibitem{BHLM2} A.~Bonami, D.~Hilhorst, E.~Logak and M.~Mimura,
{\it A free boundary problem arising in a chemotaxis model}, in
"Free Boundary Problems, Theory and Applications", M.~Niezg\'odka
and P.~Strzelecki Eds, Pitman Res. Notes in Math. Series {\bf 363}
(1996).

\bibitem{BHLM3} A.~Bonami, D.~Hilhorst, E.~Logak and M.~Mimura,
{\it Singular limit of a chemotaxis-growth model}, Advances in
Differential Equations {\bf 6} (2001), 1173--1218.

\bibitem{C1} X.~Chen,
{\it Generation and propagation of interfaces for
reaction-diffusion equations}, J.~Differential equations {\bf 96}
(1992), 116--141.

\bibitem{C2} X.~Chen, {\it Generation and propagation of interfaces
for reaction-diffusion systems}, Trans.~Amer.~Math.~Soc. {\bf 334}
(1992), 877--913.

\bibitem{CR} X.~Chen and F.~Reitich, {\it Local existence and uniqueness of solutions of the Stefan
problem with surface tension and kinetic undercooling},
J.~Math.~Anal.~Appl. {\bf 164(2)} (1992), 350--362.

\bibitem{CXY} X.~Y.~Chen, {\it Dynamics of interfaces in reaction
diffusion systems}, Hiroshima Math. J. {\bf 21} (1991), 47--83.

\bibitem{CP} S.~Childress and J.~K.~Percus, {\it Nonlinear aspects of chemotaxis},
Math. Biosci. {\bf 56} (1981), 217--237.

\bibitem{FL} R.~M.~Ford and D.~A.~Lauffenburger, {\it Analysis of chemotactic bacterial
distributions in population migration assays using a mathematical
model applicable to steep or shallow attractant gradients},
Bull.~Math.~Biol. {\bf 53} (1991), 721--749.

\bibitem{HV1} M.~A.~Herrero and J.~Velazquez, {\it Chemotactic collapse for
the Keller-Segel model}, J. Math. Biol. {\bf 35} (1996), no. 2,
177--194.

\bibitem{HV2} M.~A.~Herrero and J.~Velazquez, {\it A blow-up mechanism for a chemotaxis
model}, Ann. Scuola. Norm. Sup. Pisa Cl. Sci. {\bf 24} (1997), no.
4, 633--683.

\bibitem{JL} W.~J\"a\-ger and S.~Luckhaus, {\it On explosions of solutions to a system
of partial differential equations modelling chemotaxis}, Trans.
Amer. Math. Soc. {\bf 329} (1992), 819--824.

\bibitem{KNHM} G.~Karali, K.~Nakashima, D.~Hilhorst and H.~Matano,
{\it Singular limit of a spatially inhomogeneous Lotka-Volterra
competition-diffusion system}, in preparation.

\bibitem{KS} E.~F.~Keller and L.~A.~Segel, {\it Initiation of slime mold aggregation
viewed as an instability}, J. Theor. Biol. {\bf 26} (1970),
399--415.

\bibitem{LNT}C.~S.~Lin, W.~M.~Ni and I.~Takagi,
{\it Large amplitude stationary solutions to a chemotaxis system},
J. Differential Equations {\bf 72} (1988), 1--27.

\bibitem{MT}M.~Mimura and T.~Tsujikawa, {\it Aggregating pattern dynamics
in a  chemotaxis model including growth}, Physica A {\bf 230}
(1996), 499--543.

\bibitem{MTKU}M.~Mimura, T.~Tsujikawa, R.~Kobayashi and D.~Ueyama,
{\it Dynamics of aggregation patterns in a
chemotaxis-diffusion-growth model equation. Proceedings of the
Workshop on Principles of Pattern Formation and Morphogenesis in
Biological Systems (Kasugai, 1992/93)}, Forma {\bf 8} (1993),
179--195.

\bibitem{N} T.~Nagai, {\it Blow-up of radially symmetric solutions to
a chemotaxis pattern}, Advances in Math. Sciences and Appl. Vol.
5, no. 2 (1995), 581-601.

\bibitem{Nan} V.~Nanjundiah, {\it Chemotaxis, signal relaying and aggregation
morphology}, J. Theor. Biol. {\bf 42} (1973), 63--105.

\bibitem{NMHS} K.-I.~Nakamura, H.~Matano, D.~Hilhorst and
R.~Sch\"a\-tzle, {\it Singular limit of a reaction-diffusion
equation with a spatially inhomogeneous reaction term},
J.~Stat.~Phys. {\bf 95} (1999), 1165-1185.

\bibitem{Sch} R.~Schaaf, {\it Stationary solutions of chemotaxis systems}, Trans.
Amer. Math. Soc. {\bf 292} (1985), 531--556.

\end{thebibliography}
\end{document}